\documentclass[11pt, righteqno]{article}
\usepackage{epsfig}
\usepackage{amsfonts}
\usepackage{amsmath}
\usepackage{amssymb}
\makeindex  \textwidth 130mm

\newtheorem{theorem}{\sc Theorem}[section]
\newtheorem{proposition}{\sc Proposition}[section]
\newtheorem{corollary}{\sc Corollary}[section]

\newtheorem{remark}{\sc Remark}[section]

\def\ds {\displaystyle\mathstrut}

\providecommand{\U}[1]{\protect\rule{.1in}{.1in}}

\providecommand{\U}[1]{\protect\rule{.1in}{.1in}}

\title{Classification of
Quadratic Differential Systems on $\mathbb{R}^{3}$ Having a
Semisimple Derivation with a one-Dimensional Kernel}
\author{I. Burdujan}

\begin{document}

\title{Quadratic Differential
Systems on $\mathbb{R}^{3}$ Having a Semisimple Derivation with
one-Dimensional Kernel}
\author{I. Burdujan}

\date{}
\maketitle

\footnote{Corresponding author,

\noindent Ilie Burdujan, \emph{E-mail address}:\ $ilieburdujan@uaiasi.ro$,\ $ilieburdujan@yahoo.com$\\
address:\  3, Mihail Sadoveanu Street, 700490, Ia\c si, Romania}

\noindent\small{University of Agricultural Sciences and Veterinary
Medicine, Ia\c si, 700490, Romania}
\begin{abstract}\noindent The classification, up to a
center-affinity, of the homogeneous quadratic differential systems
defined on $\mathbb{R}^{3}$ that have at least a semisimple
derivation with one-dimensional kernel, is achieved. It is proved
that there exist 35 families of affine equivalence classes of such
systems.
\end{abstract}

\section{Introduction}

Let us recall that any homogeneous quadratic differential system
(shortly, HQDS) on an {\sc Banach} space is congenitally connected
with a commutative binary algebra. Indeed, each HQDS is defined by
means of a covariant symmetric (1,2)-tensor which, in its turn, is
the structure tensor of a commutative algebra. Consequently, the
study of any HQDS could be achieved by my means of its associated
commutative algebra. In particular, the algebra of derivations of
each HQDS is the same with derivation algebra of its associated
algebra. Of course, it is more appropriate to study the properties
of algebras having a derivation instead to make a direct study of
the corresponding HQDS. Such a study is the object of our present
paper.

Let $k$ be a field of characteristic 0 and $A(\cdot)$ be a finite
dimensional $k$-algebra. Recall that the derivation $D\in Der\ A$
is said to be \emph{semisimple} if it is diagonalisable in an
extension of $k$, i.e. there exists a basis in $A$ consisting of
eigenvectors of $D$.

In \cite{Bur5} were classified, up to an isomorphism, the real
3-dimensional commutative algebras having a semisimple nonsingular
derivation. It was shown that the existence of such a derivation
acts as a very strong constraint compelling each such algebra to
be isomorphic to one of the 4 algebras listed in \cite{Bur5}.
Accordingly, were classified the corresponding homogeneous
quadratic differential systems up to a center-affine equivalence.

The aim of this paper is to classify, up to an isomorphism, the
real 3-dimensional commutative algebras having at least a
semisimple derivation with one-dimensional kernel. The main result
is: \emph{there exist 35 families of isomorphism classes of real
3-dimensional commutative algebras having at least a semisimple
derivation with one-dimensional kernel}. For each of them are
exhibited their main properties which allow to decide on the
problem of their mutual isomorphism. Everyone of these 35 families
is either a singleton or consists of a set of algebras that have
the same lists of main properties and can be indexed by one- or
two-parameters such that two algebras in family corresponding to
different parameters are non-isomorphic.

In fact, we shall get the subalgebra lattices, the derivation
algebras and the group of automorphisms for each class of
algebras, as they are the most important invariants of binary
algebras. Especially, we are interested in finding the set $Ann\
A$ of \emph{annulator} elements, the set $\mathcal{N}(A)$ of all
\emph{nilpotent} elements and the set $\mathcal{I}(A)$ of all
\emph{idempotent} elements of algebra $A$. Further, the
corresponding homogeneous quadratic dynamical systems are
classified up to a center-affine equivalence. Recall that (see
\cite{Bur5}) the subalgebra lattice of $A(\cdot)$ allows to
identify a natural partition $\mathcal{P}_{A}$ of the ground space
$A$ which, in its turn, defines a partition of the set of all
integral curves of its associated HQDS.

\vspace{3mm}\noindent \emph{COMMENT.} Really, there exists
infinite many isomorphism classes of such algebras. In order to
identify them we need to find a partition of this set of
isomorphism classes consisting of a finite number of sets. To this
end we consider a list of "main properties" which allows to define
an equivalence, the so-called \emph{MP-equivalence}, on the set of
isomorphism classes: two isomorphism classes are MP-equivalent if
and only if their algebras have the same "main properties". The
partition associated with this equivalence has a finite number of
elements (here, 35 sets). Certainly, we can exclude in the list of
main properties the part concerning the two partitions induced by
subalgebra lattice on the ground space of analyzed algebra as well
as on the set of all integral curves of the HQDS assigned to this
algebra. We keep this part because it work like a conformity test
for automorphism group of algebra.

\section{Algebras having a semisimple derivation with one-dimensional kernel}

\setcounter{equation}{0} \vspace{3mm}\noindent Let $A(\cdot)$ be
the real 3-dimensional (nontrivial) commutative algebra associated
with a HQDS on $\mathbb{R}^{3}$. Suppose that $\widetilde{D}$ is a
nonzero semisimple derivation of $A(\cdot)$ with one-dimensional
kernel. Then, algebra $A(\cdot)$ has a semisimple derivation $D$
having the spectrum of the form $Spec\ D=(1,\omega,0)$ with
$\omega\neq 0$; this notation for spectrum is preferred because
the eigenvalues $1,\omega,0$ are not necessarily distinct each
other. Moreover, there exists a basis $\mathcal{B}=(e_{1}, e_{2},
e_{3})$ of $A$ such that
$$D(e_{1})=e_{1},\ \ \ D(e_{2})=\omega e_{2},\ \ \ D(e_{3})=0.$$
In this case, $A$ decomposes into a direct vector sum of invariant
subspaces with respect to $D$.

\vspace{3mm}In order to give the analytical expression for the
existence of a derivation $D$ of algebra $A(\cdot)$, we define -
as it is usual - the structure constants of $A$ in basis
$\mathcal{B}$ by equations
\begin{equation}\label{e21}e_{i}\cdot e_{j}=a_{ij}^{k}e_{k}.\end{equation} For convenience, we shall denote
\begin{equation}\label{e22}\begin{array}{llllll}a_{11}^{1}=a&\ \ a_{11}^{2}=b&\
\ a_{11}^{3}=c&\ \ a_{12}^{1}=k&\ \ a_{12}^{2}=m&\ \ a_{12}^{3}=n\\

a_{22}^{1}=d&\ \ a_{22}^{2}=e&\
\ a_{22}^{3}=f&\ \ a_{13}^{1}=p&\ \ a_{13}^{2}=q&\ \ a_{13}^{3}=r\\

a_{33}^{1}=g&\ \ a_{33}^{2}=h&\ \ a_{33}^{3}=j&\ \ a_{23}^{1}=s&\
\ a_{23}^{2}=t&\ \ a_{23}^{3}=v.
\end{array}\end{equation}

Then the endomorphism $D$ is a derivation for $A$ if and only if
the next conditions are fulfilled:

\begin{equation}\label{e23}\left\{\begin{array}{l}
a=c=e=f=g=h=k=m=r=v=0\\
(\omega-2)b=0\\
(1-2\omega)d=0\\
(1+\omega)n=0\\
(\omega-1)q=0\\
(1-\omega)s=0.
\end{array}\right.\end{equation}
Here $j, p, t$ range free over $\mathbb{R}$ while $b,d,n,q,s$ take
values depending on $\omega$. Equations (\ref{e23}) impose to take
into account of the natural decomposition of $\mathbb{R}^{\ast}$
as range of $\omega$ defined by means of sets:
$\{-1,\frac{1}{2},1,2\}$ and $\mathbb{R}\setminus
\{-1,0,\frac{1}{2},1,2\}$.

Consequently, we have to analyze only algebras having a semisimple
derivation $D$ with $Spec\ D$ in the following list:
$$1)\ (1,-1,0),\ \ 2)\ (1,\frac{1}{2},0),\ \ 3)\ (1,1,0),\ \ 4)\ (1,2,0),\ \ 5)\ (1,\omega,0)\ \text{with}\ \omega\notin \{-1,0,\frac{1}{2},1,2\}.$$
Since when $D\in Der\ A$ has $Spec\ D=(1,\frac{1}{2},0)$ then
derivation $D'=2D$ has $Spec\ D'=(1,2,0)$, it follows the next
result.

\begin{proposition}
If a real 3-dimensional (nontrivial) commutative algebra
$A(\cdot)$ has a semisimple derivation with one-dimensional
kernel, it has at least a derivation $D$ with $Spec\ D$ of one of
the following forms:
$$1)\ (1,-1,0),\ \  2)\ (1,1,0),\ \ 3)\ (1,2,0),\ \ 4)\ (1,\omega,0)\ with\ \omega\notin \{-1,0,\frac{1}{2},1,2\}.$$
\end{proposition}

\textbf{1) Case} $Spec\ D=(1,-1,0)$

\vspace{3mm}In basis $\mathcal{B}$ algebra $A(\cdot)$ has the next
multiplication table:
$$\begin{array}{llll}
  \textbf{Table T}\hspace{6mm} & \hspace{5mm} e_{1}^{2}=0  &\hspace{5mm} e_{2}^{2}=0&\hspace{5mm} e_{3}^{2}=je_{3} \\
   &\hspace{5mm}  e_{1}e_{2}=ne_{3}  &\hspace{5mm}  e_{1}e_{3}=pe_{1} &\hspace{5mm}  e_{2}e_{3}=te_{2}
\end{array}$$
with $j,n,p,t\in \mathbb{R}$. Algebra corresponding to $j = n
=p=t=0$ is just the null algebra that is of no interest in
general. We have to consider the case when at least one of the
parameter $j,n,p,t$ is not zero.

Further we shall consider the next two mutually exclusive cases:
$$I)\ \ jn\neq 0,\ \ \ \ \ \ II)\ \ jn=0.$$

\vspace{3mm} \emph{Case I} $jn\neq 0$

\vspace{3mm}\noindent In this case, each algebra with
multiplication table $\textbf{T}$ is isomorphic to algebra:
$$\begin{array}{llll}
  \textbf{Table T1}\hspace{6mm} & \hspace{5mm} e_{1}^{2}=0  &\hspace{5mm} e_{2}^{2}=0&\hspace{5mm} e_{3}^{2}=e_{3} \\
   &\hspace{5mm}  e_{1}e_{2}=e_{3}&\hspace{5mm}  e_{1}e_{3}=\alpha e_{1}  &\hspace{5mm}  e_{2}e_{3}=\beta e_{2}
\end{array}$$
with $\alpha,\beta\in \mathbb{R}$.

This time, $e_{3}$ is an idempotent whose left multiplication
$L_{e_{3}}$ has the spectrum $(\alpha,\beta,1)$. Consequently, as
long as $e_{3}$ is the only idempotent of such an algebra, its
eigenvalues $\alpha,\ \beta$ have to be the most important
invariants characterizing algebra.

For convenience, let us denote by $A_{1}(\alpha,\beta)$ any
algebra having the multiplication table \textbf{T1}.

\begin{proposition} The algebras $A_{1}(\alpha,\beta)$ and
$A_{1}(\beta,\alpha)$ are isomorphic.
\end{proposition}
Consequently, in the following we deal with algebras
$A_{1}(\alpha,\beta)$ with $\alpha\leq \beta$, only.

\vspace{2mm} Further, by a straightforward computation, it is
proved the next proposition.

\begin{proposition} The algebras $A_{1}(\alpha,\beta)$ (with $\alpha\leq \beta$) and
$A_{1}(\alpha_{1},\beta_{1})$ (with $\alpha_{1}\leq \beta_{1}$)
are isomorphic if and only if $\alpha=\alpha_{1}, \
\beta=\beta_{1}$.
\end{proposition}

\begin{proposition} Every algebra $A$ of type $A_{1}(\alpha,\beta)$
has:

\vspace{2mm} \ \ (i)\ $Ann\ A=\{0\}$,

\ (ii)\ $\mathcal{N}(A)=\left \{\begin{array}{llc}
\mathbb{R}e_{1}\cup \mathbb{R}e_{2}&if& \alpha^{2}+\beta^{2}\neq
0\\ \{xe_{1}-\frac{z^{2}}{2x}e_{2}+ze_{3}\ |\ x,z\in \mathbb{R},\
x\neq 0\}\cup \mathbb{R}e_{2}&if& \alpha=\beta=0
\end{array}\right.$

(iii)\ $\mathcal{I}(A)=\left \{\begin{array}{lll}\{ e_{3}\}&if&
\alpha=\beta=0\\
\{e_{3}\}&if&
\alpha\notin\{0,\frac{1}{2}\},\ \beta=0\\
 \{xe_{1}+e_{3}\ | \ x\in\mathbb{R}\}&if&
\alpha=\frac{1}{2},\ \beta\neq \frac{1}{2}\\
\{e_{3}\}&if&
\alpha=0,\ \beta\notin\{0,\frac{1}{2}\}\\
 \{ye_{2}+e_{3}\ | \ y\in\mathbb{R}\}&if&\alpha\neq \frac{1}{2},\ \beta=\frac{1}{2}\\

\{xe_{1}+\frac{2\alpha-1}{8\alpha^{2}}\frac{1}{x}e_{2}+\frac{1}{2\alpha}e_{3}\
| \ x\in\mathbb{R}^{\ast}\}\cup\{e_{3}\}&if&
\alpha=\beta\notin\{0,\ \frac{1}{2}\}\\

 \{xe_{1}+e_{3}\ |\ x\in\mathbb{R}\}\cup\{ye_{2}+e_{3}\ |\ y\in \mathbb{R}\}&if&
\alpha=\beta=\frac{1}{2}\\
\{e_{3}\}&if& \alpha\neq \beta,\ \alpha,\
\beta\notin\{0,\frac{1}{2}\}.
\end{array}\right.$
\end{proposition}

We start the study of each kind of algebra exhibited in previous
Proposition, taking into account that algebras $A_{1}(\alpha,0)$
and $A_{1}(0,\alpha)$ are isomorphic.

\vspace{3mm}\textbf{A1)} \emph{Properties of algebra
$A=A_{1}(0,0)$}

\vspace{3mm}$\bullet$ $Ann\ A=\{0\},\
\mathcal{N}(A)=\{xe_{1}-\frac{z^{2}}{2x}e_{2}+ze_{3}\ |\ x,z\in
\mathbb{R},\ z\neq 0\}\cup \mathbb{R}e_{2},\
\mathcal{I}(A)=\{e_{3}\}$,

$\bullet$ 1-dimensional subalgebras: $\mathbb{R}u$ for
$u\in\mathcal{N}(A)\cup \mathcal{I}(A)$,

$\bullet$ 2-dimensional subalgebras:
$Span_{\mathbb{R}}\{e_{3},ae_{1}+be_{2}\}\ (a^{2}+b^{2}\neq 0)$,

$\bullet$ ideals: $\mathbb{R}e_{3},\
Span_{\mathbb{R}}\{e_{3},ae_{1}+be_{2}\}\ (a^{2}+b^{2}\neq 0)$,

$\bullet$ $A^{2}=\mathbb{R}e_{3}$; $A/A^{2}$ is the null
2-dimensional algebra,

$\bullet$ $Der\ A= \mathbb{R}D$,

$\bullet$ $Aut\ A=H\cup JH$ where $$H= \left\{
\left[\begin{array}{lll}x&0&0\\0&x^{-1}&0\\0&0&1\end{array}\right]\
|\ x\in \mathbb{R}^{\ast}\right\},\
J=\left[\begin{array}{lll}0&1&0\\1&0&0\\0&0&1\end{array}\right] \
\ \text{and}\ \ J^{2}=id;$$ $H$ is a normal divisor of $Aut\ A$
and $Aut\ A/H\cong \mathbb{Z}_{2}$ (in fact, $Aut\ A\cong
\mathbb{R}^{\ast}\times \{-1,1\}$ where $\mathbb{R}^{\ast}$
 denotes the multiplicative group of
nonzero real numbers $\mathbb{R}^{\ast}(\cdot)$); $H$ is a normal
divisor of $Aut\ A$,

$\bullet$ the partition $\mathcal{P}_{A}$ of $\mathbb{R}^{3}$,
defined by the lattice of subalgebras of $A$, consists of:

\vspace{1mm} \hspace{5mm} $\diamond$  the singletons covering the
lines $\mathbb{R}u$ for $u\in \mathcal{N}(A)$ (i.e. covering the
cone $2x^{1}x^{2}+(x^{3})^{2}=0$),

\hspace{5mm} $\diamond$  the half-lines of axis $Ox^{3}$ delimited
by $O$,

\hspace{5mm} $\diamond$  the connected components of each plane
passing through $Ox^{3}$, delimited by axis $Ox^{3}$ and the
generatrices of cone  $2x^{1}x^{2}+(x^{3})^{2}=0$ (whenever it is
the case, i.e. $x^{1}x^{2}<0$),

\vspace{1mm}$\bullet$ the partition $\mathcal{P}_{A}$ of $A$
induces a partition on the set of integral curves of the
associated homogeneous quadratic differential system (HQDS)
consisting of:

\vspace{1mm} \hspace{5mm} $\diamond$  the singletons consisting of
singular solutions that cover the cone
$2x^{1}x^{2}+(x^{3})^{2}=0$,

\hspace{5mm} $\diamond$  the families of ray-solutions contained
in each half-line of $Ox^{3}$ delimited by $O$,

\hspace{5mm} $\diamond$  the integral curves contained in the
connected components of each plane passing through $Ox^{3}$,
delimited by axis $Ox^{3}$ and the generatrices of cone
$2x^{1}x^{2}+(x^{3})^{2}=0$ (whenever it is the case, i.e.
$x^{1}x^{2}<0$).

\vspace{3mm}Since each nonsingular integral curves lies into a
plane passing through $Ox^{3}$ it has a null torsion tensor.
Moreover, $A/A^{2}$ is a null algebra what implies that the
curvature tensor of nonsingular integral curves vanishes too, so
that each non-singular integral curves lies on a line parallel to
$Ox^{3}$.

\vspace{3mm}\textbf{A2)} \emph{Properties of algebra
$A=A_{1}(0,\frac{1}{2})$}

\vspace{3mm}$\bullet$ $Ann\ A=\{0\},\
\mathcal{N}(A)=\mathbb{R}e_{1}\cup \mathbb{R}e_{2},\
\mathcal{I}(A)=\{ye_{2}+e_{3}\ |\ y\in \mathbb{R}\}$,

$\bullet$ 1-dimensional subalgebras: $\mathbb{R}u$ for $u\in
\mathcal{N}(A)\cup \mathcal{I}(A)$,

$\bullet$ 2-dimensional subalgebras:
$Span_{\mathbb{R}}\{e_{1},e_{3}\},\
Span_{\mathbb{R}}\{e_{2},e_{3}\}$,

$\bullet$ ideals: $Span_{\mathbb{R}}\{e_{2},e_{3}\}$,

$\bullet$ $A^{2}=Span_{\mathbb{R}}\{e_{2},e_{3}\}$; $A/A^{2}$ is
the null 1-dimensional algebra,

$\bullet$ $Der\ A= \mathbb{R}D$,

$\bullet$ $Aut\ A= \left\{
\left[\begin{array}{lll}x&0&0\\0&x^{-1}&0\\0&0&1\end{array}\right]\
|\ x\in \mathbb{R}^{\ast}\right\}\cong\mathbb{R}^{\ast}(\cdot)$.

$\bullet$ the partition $\mathcal{P}_{A}$ of $\mathbb{R}^{3}$,
defined by the lattice of subalgebras of $A$, consists of:

\vspace{1mm} \hspace{5mm} $\diamond$  the singletons covering the
axes $\mathbb{R}e_{1}$ and $\mathbb{R}e_{2}$,

\hspace{5mm} $\diamond$  the half-lines delimited by $O$ on lines
$\mathbb{R}u$ for $u\in \mathcal{I}(A)$ (i.e. these lines cover
the plane $x^{2}Ox^{3}$ less $x^{2}$-axis),

\hspace{5mm} $\diamond$ the half-spaces $x^{1}<0$ and $x^{1}>0$
(delimited by plane $x^{2}Ox^{3}$) less the points of axis
$Ox^{1}$,

\vspace{1mm}$\bullet$ the partition $\mathcal{P}_{A}$ of $A$
induces a partition on the set of integral curves of the
associated homogeneous quadratic differential system (HQDS)
consisting of:

\vspace{1mm} \hspace{5mm} $\diamond$  the singletons consisting of
singular solutions that cover the axes $\mathbb{R}e_{1}$ and
$\mathbb{R}e_{2}$,

\hspace{5mm} $\diamond$  the families of ray-solutions contained
in each half-line delimited by $O$ on each line $\mathbb{R}u$ for
$u\in \mathcal{I}(A)$,

\hspace{5mm} $\diamond$  the integral curves contained in the
half-spaces $x^{1}<0$ and $x^{1}>0$ (delimited by plane
$x^{2}Ox^{3}$) less $x^{1}$-axis.

Let us remark that each nonsingular integral curve has a null
torsion tensor (indeed, $A^{2}$ is a proper ideal of $A$).

\begin{remark} All idempotents have the same spectrum
$(0,\frac{1}{2},1)$. Then $A$ has the decomposition $A(0)\oplus
A(\frac{1}{2})\oplus A(1)$ consisting of the eigenspaces
corresponding to any idempotent used in basis.\end{remark}

\vspace{2mm}\textbf{A3)} \emph{Properties of algebras $A$ of type
$A_{1}(\frac{1}{2},\frac{1}{2})$}

\vspace{2mm}$\bullet$ $Ann\ A=\{0\},\
\mathcal{N}(A)=\mathbb{R}e_{1}\cup \mathbb{R}e_{2},\
\mathcal{I}(A)=\{xe_{1}+e_{3}\ |\ x\in \mathbb{R}\}\cup
\{ye_{2}+e_{3}\ |\ y\in \mathbb{R}^{\ast}\}$,

$\bullet$ 1-dimensional subalgebras: $\mathbb{R}u$ for
$u\in\mathcal{N}(A)\cup \mathcal{I}(A)$,

$\bullet$ 2-dimensional subalgebras: $Span_{\mathbb{R}}\{e_{3},
pe_{1}+qe_{2}\}$ with $p,q\in \mathbb{R}$,

$\bullet$ ideals: none,

$\bullet$ $A^{2}=A$,

$\bullet$ $Der\ A= \mathbb{R}D$,

$\bullet$ $Aut\ A=H\cup JH$ where $$H= \left\{
\left[\begin{array}{lll}x&0&0\\0&x^{-1}&0\\0&0&1\end{array}\right]\
|\ x\in \mathbb{R}^{\ast}\right\},\
J=\left[\begin{array}{lll}0&1&0\\1&0&0\\0&0&1\end{array}\right] \
\ \text{and}\ \ J^{2}=id,$$ (i.e. $Aut\ A\cong
\mathbb{R}^{\ast}\times \{-1,1\}$) and $H$ is a normal divisor of
$Aut\ A$,

$\bullet$ the partition $\mathcal{P}_{A}$ of $\mathbb{R}^{3}$,
defined by the lattice of subalgebras of $A$, consists of:

\vspace{1mm} \hspace{5mm} $\diamond$  the singletons covering the
axes $Ox^{1}$ and $Ox^{2}$,

\hspace{5mm} $\diamond$  the half-lines of $x^{1}Ox^{3}$ passing
through $O$, delimited by $O$, less the axis $Ox^{1}$,

\hspace{5mm} $\diamond$  the half-lines of $x^{2}Ox^{3}$ passing
through $O$, delimited by $O$, less the axes $Ox^{2}$ and
$Ox^{3}$,

\hspace{5mm} $\diamond$  the half-planes delimited by axis
$Ox^{3}$ on each plan containing axis $Ox^{3}$, except the planes
$x^{1}Ox^{3}$ and $x^{2}Ox^{3}$,

$\bullet$ the partition $\mathcal{P}_{A}$ of $A$ induces a
partition on the set of integral curves of the associated
homogeneous quadratic differential system (HQDS) consisting of:

\vspace{1mm} \hspace{5mm} $\diamond$  the singletons consisting of
singular solutions that cover the axes $Ox^{1}$ and $Ox^{2}$,

\hspace{5mm} $\diamond$  the families of ray-solutions contained
in each half-line of $x^{1}Ox^{3}$ passing through $O$ and
delimited by $O$, less the axis $Ox^{1}$,

\hspace{5mm} $\diamond$  the families of ray-solutions contained
in each half-line of $x^{2}Ox^{3}$ passing through $O$ and
delimited by $O$, less the axes $Ox^{2}$ and $Ox^{3}$,

\hspace{5mm} $\diamond$  the integral curves contained in each
half-plane delimited by axis $Ox^{3}$ on each plan containing axis
$Ox^{3}$ without planes $x^{1}Ox^{3}$ and $x^{2}Ox^{3}$.

\vspace{3mm}Note that each nonsingular solution is a torsion-free
curve. Moreover, all idempotents have the same spectrum
$(\frac{1}{2},\frac{1}{2},1)$.

\vspace{3mm}\textbf{A4)} \emph{Properties of algebra $A$ of type
$A_{1}(\alpha,0)$ with $\alpha\notin\{0,\frac{1}{2}\}$}

\vspace{2mm}$\bullet$ the algebras $A_{1}(\alpha,0)$ and
$A_{1}(\alpha',0)$ are isomorphic if and only if $\alpha=\alpha'$,

$\bullet$ $Ann\ A=\{0\},\ \mathcal{N}(A)=\mathbb{R}e_{1}\cup
\mathbb{R}e_{2},\ \mathcal{I}(A)=\{e_{3}\}$,

$\bullet$ 1-dimensional subalgebras: $\mathbb{R}e_{1},\
\mathbb{R}e_{2},\ \mathbb{R}e_{3}$,

$\bullet$ 2-dimensional subalgebras:
$Span_{\mathbb{R}}\{e_{1},e_{3}\},\
Span_{\mathbb{R}}\{e_{2},e_{3}\}$,

$\bullet$ ideals: $Span_{\mathbb{R}}\{e_{1},e_{3}\}$,

$\bullet$ $A^{2}=Span_{\mathbb{R}}\{e_{1},e_{3}\}$; $A/A^{2}$ is
the null 1-dimensional algebra,

$\bullet$ $Der\ A= \mathbb{R}D$,

$\bullet$ $Aut\ A= \left\{
\left[\begin{array}{lll}x&0&0\\0&x^{-1}&0\\0&0&1\end{array}\right]\
|\ x\in \mathbb{R}^{\ast}\right\}\cong \mathbb{R}^{\ast}(\cdot)$.

$\bullet$ the partition $\mathcal{P}_{A}$ of $\mathbb{R}^{3}$,
defined by the lattice of subalgebras of $A$, consists of:

\vspace{1mm} \hspace{5mm} $\diamond$  the singletons covering the
axes $\mathbb{R}e_{1}$ and $\mathbb{R}e_{2}$,

\hspace{5mm} $\diamond$  the half-lines of axis $Ox^{3}$ delimited
by $O$,

\hspace{5mm} $\diamond$ the quarters of plane $x^{1}Ox^{3}$
delimited by axes $Ox^{1}$ and $Ox^{3}$,

\hspace{5mm} $\diamond$ the quarters of plane $x^{2}Ox^{3}$
delimited by axes $Ox^{2}$ and $Ox^{3}$,

\hspace{5mm} $\diamond$ the quarters of space delimited by planes
$x^{1}Ox^{3}$ and $Ox^{2}Ox^{3}$,

\vspace{1mm}$\bullet$ the partition $\mathcal{P}_{A}$ of $A$
induces a partition on the set of integral curves of the
associated homogeneous quadratic differential system (HQDS)
consisting of:

\vspace{1mm} \hspace{5mm} $\diamond$  the singletons consisting of
singular solutions that cover the axes $\mathbb{R}e_{1}$ and
$\mathbb{R}e_{2}$,

\hspace{5mm} $\diamond$  the families of ray-solutions contained
in each half-line of $Ox^{3}$ delimited by $O$,

\hspace{5mm} $\diamond$  the integral curves contained in the
quarters of plane $x^{1}Ox^{3}$ delimited by axes $Ox^{1}$ and
$Ox^{3}$,

\hspace{5mm} $\diamond$  the integral curves contained in the
quarters of plane $x^{2}Ox^{3}$ delimited by axes $Ox^{2}$ and
$Ox^{3}$,

\hspace{5mm} $\diamond$  the integral curves contained in the
quarters of space delimited by planes $x^{1}Ox^{3}$ and
$Ox^{2}Ox^{3}$.

Note that each nonsingular integral curve has a null torsion
tensor.

\vspace{3mm}\noindent\textbf{Note.} In order to bring together the
classes $A_{1}(\alpha,0)$ with $\alpha<0$ and $A_{1}(0,\beta)$
with $\beta>0$ we have ignored, in case \textbf{A4)}, the
convention to consider increasing parameters in algebras of type
$A_{1}(a,b)$. Such a behavior will be used whenever is possible,
in order to save the space.

\vspace{5mm} \textbf{A5)}  \emph{Properties of algebras $A$ of
type $A_{1}(\alpha,\frac{1}{2})$ with}
$\alpha\notin\{0,\frac{1}{2}\}$

\vspace{3mm}$\bullet$ algebras $A_{1}(\alpha,\frac{1}{2})\
(\alpha\notin\{0,\frac{1}{2}\})$ and $A_{1}(\alpha',\frac{1}{2})\
(\alpha'\notin\{0,\frac{1}{2}\})$ are isomorphic if and only if
$\alpha=\alpha'$,

$\bullet$ $Ann\ A=\{0\},\ \mathcal{N}(A)=\mathbb{R}e_{1}\cup
\mathbb{R}e_{2},\ \mathcal{I}(A)=\{ye_{2}+e_{3}\ |\ y\in
\mathbb{R}\}$,

$\bullet$ 1-dimensional subalgebras: $\mathbb{R}u$ for
$u\in\mathcal{N}(A)\cup \mathcal{I}(A)$,

$\bullet$ 2-dimensional subalgebras:
$Span_{\mathbb{R}}\{e_{1},e_{3}\}$,
$Span_{\mathbb{R}}\{e_{2},e_{3}\}$,

$\bullet$ ideals: none,

$\bullet$ $A^{2}=A$,

$\bullet$ $Der\ A= \mathbb{R}D$,

$\bullet$ $Aut\ A= \left\{
\left[\begin{array}{lll}x&0&0\\0&x^{-1}&0\\0&0&1\end{array}\right]\
|\ x\in \mathbb{R}^{\ast}\right\}\cong \mathbb{R}^{\ast}(\cdot)$,

$\bullet$ the partition $\mathcal{P}_{A}$ of $\mathbb{R}^{3}$,
defined by the lattice of subalgebras of $A$, consists of:

\vspace{1mm} \hspace{5mm} $\diamond$  the singletons covering the
axes $Ox^{1}$ and $Ox^{2}$,

\hspace{5mm} $\diamond$  the half-lines delimited by $O$ on lines
of $x^{2}Ox^{3}$ passing through $O$, less the axis $Ox^{2}$,

\hspace{5mm} $\diamond$ the quarters of plane $x^{1}Ox^{3}$,
delimited by axis $Ox^{1}$ and $Ox^{3}$,

\hspace{5mm} $\diamond$ the quarter-spaces delimited by planes
$x^{1}Ox^{3}$ and $x^{2}Ox^{3}$,

\vspace{1mm}$\bullet$ the partition $\mathcal{P}_{A}$ of $A$
induces a partition on the set of integral curves of the
associated homogeneous quadratic differential system (HQDS)
consisting of:

\vspace{1mm} \hspace{5mm} $\diamond$  the singletons consisting of
singular solutions that cover the axes $Ox^{1}$ and $Ox^{2}$,

\hspace{5mm} $\diamond$  the families of ray-solutions contained
in each half-line delimited by $O$ on each line in $x^{2}Ox^{3}$
passing through $O$, less the axis $Ox^{2}$,

\hspace{5mm} $\diamond$  the integral curves contained in each
quarter of plane $x^{1}Ox^{3}$, delimited by axes $Ox^{1}$ and
$Ox^{3}$,

\hspace{5mm} $\diamond$ the integral curves contained in each
quarter-space delimited by planes $x^{1}Ox^{3}$ and $x^{2}Ox^{3}$.

\vspace{3mm}Moreover, all idempotents have the same spectrum
$(\alpha, \frac{1}{2},1)$.

\vspace{3mm}\noindent\textbf{Note.} In order to bring together the
classes $A_{1}(\alpha,\frac{1}{2})$ and $A_{1}(\frac{1}{2},\beta)$
we have ignored, in case \textbf{A5}, the convention to consider
increasing parameters in algebras of type $A_{1}(a,b)$.

\vspace{3mm} \textbf{A6)}  \emph{Properties of algebras $A$ of
type $A_{1}(\alpha,\alpha)$ with} $\alpha\notin \{0,\frac{1}{2}\}$

\vspace{3mm}$\bullet$ algebras $A_{1}(\alpha,\alpha)\
(\alpha\notin \{0,\frac{1}{2}\})$ and $A_{1}(\alpha',\alpha')\
(\alpha'\notin \{0,\frac{1}{2}\})$ are isomorphic if and only if
$\alpha=\alpha'$,

$\bullet$ $Ann\ A=\{0\},\ \mathcal{N}(A)=\mathbb{R}e_{1}\cup
\mathbb{R}e_{2},\ \mathcal{I}(A)=\{e_{3}\}\cup
\{xe_{1}+\frac{2\alpha-1}{8\alpha^{2}x}e_{2}+\frac{1}{2\alpha}e_{3}\
|\ x\in \mathbb{R}^{\ast}\}$ (each idempotent lies on the cone
$2x^{1}x^{2}+(1-2\alpha)(x^{3})^{2}=0$),

$\bullet$ 1-dimensional subalgebras: $\mathbb{R}u$ for
$u\in\mathcal{N}(A)\cup \mathcal{I}(A)$,

$\bullet$ 2-dimensional subalgebras:
$Span_{\mathbb{R}}\{e_{3},ae_{1}+be_{2}\}$ for $a^{2}+b^{2}\neq
0$,

$\bullet$ ideals: none,

$\bullet$ $A^{2}=A$,

$\bullet$ $Der\ A= \mathbb{R}D$,

$\bullet$ $Aut\ A=H\cup JH$ where $$ H= \left\{
\left[\begin{array}{lll}x&0&0\\0&x^{-1}&0\\0&0&1\end{array}\right]\
|\ x\in \mathbb{R}^{\ast}\right\},\ \
J=\left[\begin{array}{lll}0&1&0\\1&0&0\\0&0&1\end{array}\right] \
\ \text{and}\ \ J^{2}=id,$$ i.e. $Aut\ A\cong
\mathbb{R}^{\ast}\times \{-1,1\}$ and $H$ is a normal divisor of
$Aut\ A$,

$\bullet$ the partition $\mathcal{P}_{A}$ of $\mathbb{R}^{3}$,
defined by the lattice of subalgebras of $A$, consists of:

\vspace{1mm} \hspace{5mm} $\diamond$  the singletons covering the
axes $Ox^{1}$ and $Ox^{2}$,

\hspace{5mm} $\diamond$  the half-lines of axis $Ox^{3}$ delimited
by $O$,

\hspace{5mm} $\diamond$ the half-lines delimited by $O$ on each
generatrix of cone $2x^{1}x^{2}+(1-2\alpha)(x^{3})^{2}=0$,

\hspace{5mm} $\diamond$ the quarters of plane $x^{1}Ox^{3}$
delimited by axes $Ox^{1}$ and  $Ox^{3}$,

\hspace{5mm} $\diamond$ the quarters of plane $x^{2}Ox^{3}$
delimited by axes $Ox^{2}$ and  $Ox^{3}$,

\hspace{5mm} $\diamond$ the connected components of planes passing
through $Ox^{3}$, without planes $x^{1}Ox^{3}$ and $x^{2}Ox^{3}$,
delimited by axis $Ox^{3}$ and the cone
$2x^{1}x^{2}+(1-2\alpha)(x^{3})^{2}=0$ (whenever it is the case,
i.e. if $(1-2\alpha)x^{1}x^{2}<0)$,

$\bullet$ the partition $\mathcal{P}_{A}$ of $A$ induces a
partition on the set of integral curves of the associated
homogeneous quadratic differential system (HQDS) consisting of:

\vspace{1mm} \hspace{5mm} $\diamond$  the singletons consisting of
singular solutions that cover the axes $Ox^{1}$ and $Ox^{2}$,

\hspace{5mm} $\diamond$  the families of ray solutions contained
in the half-lines of axis $Ox^{3}$ delimited by $O$,

\hspace{5mm} $\diamond$  the families of ray solutions contained
in the half-lines delimited by $O$ on each generatrix of cone
$2x^{1}x^{2}+(1-2\alpha)(x^{3})^{2}=0$,

\hspace{5mm} $\diamond$ the integral curves contained in each
quarter of plane $x^{1}Ox^{3}$ delimited by axes $Ox^{1}$ and
 $Ox^{2}$,

\hspace{5mm} $\diamond$  the integral curves contained in each
quarter of plane $x^{2}Ox^{3}$ delimited by axis $Ox^{2}$,

\hspace{5mm} $\diamond$ the integral curves contained in the
connected components of planes passing through $Ox^{3}$, without
planes $x^{1}Ox^{3}$ and $x^{2}Ox^{3}$, delimited by axis $Ox^{3}$
and the cone $2x^{1}x^{2}+(1-2\alpha)(x^{3})^{2}=0$ (whenever it
is the case, i.e. if $(1-2\alpha)x^{1}x^{2}<0)$.

\vspace{3mm}Note that, each non-singular integral curve has
torsion zero. Moreover, idempotent $e_{3}$ has the spectrum
$\{\alpha,\alpha,1\}$ while all other idempotents have the same
spectrum $\{\frac{1}{2},1,\frac{1-\alpha}{2\alpha}\}$ what could
be connected with the presence of two connected components of
$Aut\ A$.

\vspace{5mm} \textbf{A7)}  \emph{Properties of algebras $A$ of
type $A_{1}(\alpha,\beta)$ with} $\alpha<\beta$ and
$\alpha,\beta\notin \{0, \frac{1}{2}\}$

\vspace{3mm}$\bullet$ algebras $A_{1}(\alpha,\beta)\
(\alpha<\beta)$ and $A_{1}(\alpha',\beta')\ (\alpha'<\beta')$ are
isomorphic if and only if $\alpha=\alpha',\ \ \beta=\beta'$,

$\bullet$ $Ann\ A=\{0\},\ \mathcal{N}(A)=\mathbb{R}e_{1}\cup
\mathbb{R}e_{2},\ \mathcal{I}(A)=\{e_{3}\}$,

$\bullet$ 1-dimensional subalgebras: $\mathbb{R}e_{1},\
\mathbb{R}e_{2},\ \mathbb{R}e_{3}$,

$\bullet$ 2-dimensional subalgebras: $Span_{\mathbb{R}}\{e_{1},
e_{3}\}$, $Span_{\mathbb{R}}\{e_{2}, e_{3}\}$,

$\bullet$ ideals: none,

$\bullet$ $A^{2}=A$,

$\bullet$ $Der\ A= \mathbb{R}D$,

$\bullet$ $Aut\ A= \left\{
\left[\begin{array}{lll}x&0&0\\0&x^{-1}&0\\0&0&1\end{array}\right]\
|\ x\in \mathbb{R}^{\ast}\right\}\cong \mathbb{R}^{\ast}(\cdot)$.

$\bullet$ the partition $\mathcal{P}_{A}$ of $\mathbb{R}^{3}$,
defined by the lattice of subalgebras of $A$, consists of:

\vspace{1mm} \hspace{5mm} $\diamond$  the singletons covering the
axes $Ox^{1}$ and $Ox^{2}$,

\hspace{5mm} $\diamond$ the half-axes delimited by
 $O$ on $Ox^{3}$,

\hspace{5mm} $\diamond$  the quarters delimited by axes $Ox^{1}$
and $Ox^{3}$ on plane $x^{1}Ox^{3}$,

\hspace{5mm} $\diamond$ the the quarters delimited by axes
$Ox^{2}$ and $Ox^{3}$ on plane $x^{2}Ox^{3}$,

\hspace{5mm} $\diamond$ the quarter-spaces delimited by planes
$x^{1}Ox^{3}$ and $x^{2}Ox^{3}$,

\vspace{1mm}$\bullet$ the partition $\mathcal{P}_{A}$ of $A$
induces a partition on the set of integral curves of the
associated homogeneous quadratic differential system (HQDS)
consisting of:

\vspace{1mm} \hspace{5mm} $\diamond$  the singletons consisting of
singular solutions that cover the axes $Ox^{1}$ and $Ox^{2}$,

\hspace{5mm} $\diamond$  the ray solutions lying on semi-axes of
$Ox^{3}$, delimited by $O$,

\hspace{5mm} $\diamond$  the integral curves contained in each
quarters of plane of $x^{1}Ox^{3}$, delimited by axes $Ox^{1}$ and
$Ox^{3}$,

\hspace{5mm} $\diamond$  the integral curves contained in each
quarter of plane  of $x^{2}Ox^{3}$, delimited by axis $Ox^{2}$ and
$Ox^{3}$,

\hspace{5mm} $\diamond$ the integral curves contained in each
quarter-space delimited by planes $x^{1}Ox^{3}$ and $x^{2}Ox^{3}$.

\vspace{3mm} By comparing the lists of properties of algebras in
classes \textbf{A1-A7} and taking into account of remarks
concerning the spectrum of their idempotents it follows the next
result.
\begin{theorem}\label{t1} Each algebra of type $\textbf{Ai}$ is
not isomorphic to any algebra of type $\textbf{Aj}$ for
$i,j\in\{1,2,...,7\}$ and $i\neq j$.\end{theorem}

\vspace{3mm} \emph{Case II} $jn\ =\ 0$

\vspace{3mm} This case decomposes in:

\vspace{2mm} \ \ $(i)$\ \ $j\neq 0,\ \ n=0,$

\ $(ii)$\ \ $j =0,\ \ n\neq 0,$

$(iii)$\ \ $j =0,\ \ n= 0.$

\vspace{3mm} \emph{Subcase (i)\ \ $j\neq 0$,\ \ n=0}

\vspace{3mm} There exists a basis such that the multiplication
table of algebra becomes

$$\begin{array}{llll}
  \textbf{Table T2}\hspace{6mm} & \hspace{5mm} e_{1}^{2}=0  &\hspace{5mm} e_{2}^{2}=0&\hspace{5mm} e_{3}^{2}=e_{3} \\
   &\hspace{5mm}  e_{1}e_{2}=0&\hspace{5mm}  e_{1}e_{3}=\alpha e_{1}  &\hspace{5mm}  e_{2}e_{3}=\beta e_{2}
\end{array}$$
with $\alpha,\ \beta\in \mathbb{R}$. Let us denote by
$A_{2}(\alpha, \beta)$ any algebra of type $\textbf{T2}$.

\begin{proposition}
The algebras $A_{2}(\alpha, \beta)$ and $A_{2}( \beta, \alpha)$
are isomorphic.
\end{proposition}
Consequently, in the following we restrict our interest to the
case $\alpha\leq \beta$.

\begin{proposition}
The algebras $A_{2}(\alpha, \beta)$ with $\alpha\leq \beta$ and
$A_{2}(\alpha', \beta')$ with $\alpha'\leq \beta'$ are isomorphic
if and only if $\alpha =\alpha'$ and $\beta =\beta'$.
\end{proposition}

The next result holds true for any algebra $A$ of type
$A_{2}(\alpha,\beta)$.

\begin{proposition} Every algebra $A$ of type $A_{2}(\alpha,\beta)$
has:
$$Ann\ A=\left\{ \begin{array}{lcll}(i_{1})&Span_{\mathbb{R}}\{e_{1}, e_{2}\} &if& \alpha=\beta=0\\
(i_{2})&\mathbb{R}e_{1}&if& \alpha=0,\ \beta\neq 0\\ (i_{3})&\mathbb{R}e_{2}&if& \alpha\neq 0,\ \beta= 0\\
(i_{4})&\{0\}&if& \alpha\beta\neq 0.
\end{array}\right.$$

$$\mathcal{N}(A)=Span_{\mathbb{R}}\{e_{1}, e_{2}\}$$

$$\mathcal{I}(A)=\left\{\begin{array}{lll}\{e_{3}\}&if&\alpha\neq \frac{1}{2},\ \beta\neq \frac{1}{2}\\
\{xe_{1}+e_{3}\ |\ x\in \mathbb{R}\}&if&\alpha= \frac{1}{2},\ \beta\neq \frac{1}{2}\\
\{ye_{2}+e_{3}\ |\ x\in \mathbb{R}\}&if&\alpha\neq \frac{1}{2},\ \beta= \frac{1}{2}\\
\{xe_{1}+ye_{2}+e_{3}\ |\ x,y\in \mathbb{R}\}&if&\alpha=\beta=
\frac{1}{2}
\end{array}\right.$$
\end{proposition}

\vspace{5mm}\textbf{A8)} \emph{Properties of algebras
$A=A_{2}(0,0)$}

\vspace{2mm}$\bullet$ $Ann\ A=Span_{\mathbb{R}}\{e_{1}, e_{2}\},\
\mathcal{N}(A)=Span_{\mathbb{R}}\{e_{1}, e_{2}\},\
\mathcal{I}(A)=\{e_{3}\}$,

$\bullet$ 1-dimensional subalgebras: $\mathbb{R}u$ for $u\in
\mathcal{N}(A)\cup \mathcal{I}(A)$,

$\bullet$ 2-dimensional subalgebras: $Span_{\mathbb{R}}\{e_{3},
pe_{1}+ qe_{2}\}$, $Span_{\mathbb{R}}\{e_{1}, e_{2}\}$,

$\bullet$ ideals: $\mathbb{R}e_{3},\ \mathbb{R}(pe_{1}+ qe_{2}),\
Span_{\mathbb{R}}\{e_{1}, e_{2}\},\ Span_{\mathbb{R}}\{e_{3},
pe_{1}+ qe_{2}\}\ (p^{2}+q^{2}\neq 0),$

$\bullet$ $A^{2}=\mathbb{R}e_{3}$; $A/A^{2}$ is the null
2-dimensional algebra and $A=Span_{\mathbb{R}}\{e_{1},
e_{2}\}\oplus \mathbb{R}e_{3}$ is a direct sum of two ideals (this
is like a {\sc Wedderburn-Artin} decomposition by means of a
maximal nilpotent ideal and a complementary subalgebra); $A$ is an
associative algebra,

$\bullet$ $Der\ A= \left\{\left[\begin{array}{lll} x&u&0\\
z&u&0\\0&0&0
\end{array} \right]\ |\
 x,y,z,u \in \mathbb{R}\right\}$ (i.e. $Der\ A\cong g\ell(2,\mathbb{R})$),

$\bullet$ $Aut\ A= \left\{
\left[\begin{array}{lll}x&y&0\\z&u&0\\0&0&1\end{array}\right]\ |\
x,y,z,u\in \mathbb{R},\ xu-yz\neq 0\right\}$ (i.e $Aut\ A\cong
GL(2,\mathbb{R})$),

$\bullet$ the partition $\mathcal{P}_{A}$ of $\mathbb{R}^{3}$,
defined by the lattice of subalgebras of $A$, consists of:

\vspace{1mm} \hspace{5mm} $\diamond$  the singletons covering the
plane $x^{1}Ox^{2}$,

\hspace{5mm} $\diamond$  the half-axes of $Ox^{3}$, delimited by
 $O$,

\hspace{5mm} $\diamond$  the quarters of planes delimited by axis
$Ox^{3}$ and plane $x^{1}Ox^{2}$ on each plane containing axis
$Ox^{3}$,

\vspace{1mm}$\bullet$ the partition $\mathcal{P}_{A}$ of $A$
induces a partition on the set of integral curves of the
associated homogeneous quadratic differential system (HQDS)
consisting of:

\vspace{1mm} \hspace{5mm} $\diamond$ the singletons consisting of
singular solutions that cover the plane $x^{1}Ox^{2}$,

\hspace{5mm} $\diamond$  the families of ray solutions lying on
semi-axes of $Ox^{3}$ delimited by $O$,

\hspace{5mm} $\diamond$  the integral curves contained in each
quarters of each plane passing through axis $Ox^{1}$, delimited by
axis $Ox^{3}$ and plane $x^{1}Ox^{2}$.

\vspace{3mm}Note that each integral curve lies on a half-line
parallel to $Ox^{3}$ delimited by plane $x^{1}Ox^{2}$. It means
that each integral curve has both curvature and torsion tensors
zero.

\vspace{5mm}\textbf{A9)} \emph{Properties of algebras
$A=A_{2}(0,\frac{1}{2})$}

\vspace{3mm}$\bullet$ $Ann\ A=\mathbb{R}e_{1},\ \mathcal{N}(A)=
Span_{\mathbb{R}}\{e_{1}, e_{2}\},\ \mathcal{I}(A)=\{ye_{2}+e_{3}\
|\ y\in \mathbb{R}\}$,

$\bullet$ 1-dimensional subalgebras: $\mathbb{R}u$ for $u\in
\mathcal{N}(A)\cup \mathcal{I}(A)$,

$\bullet$ 2-dimensional subalgebras: $Span_{\mathbb{R}}\{e_{1},
be_{2}+ce_{3}\}\ (b^{2}+c^{2}\neq 0)$, $Span_{\mathbb{R}}\{e_{2},
e_{3}\}$,

$\bullet$ ideals: $\mathbb{R}e_{1}$, $Span_{\mathbb{R}}\{e_{2},
e_{3}\}$,

$\bullet$ $A^{2}=Span_{\mathbb{R}}\{e_{2}, e_{3}\}$; $A/A^{2}$ is
the null 1-dimensional algebra; $A$ is a vector direct sum of
ideals,

$\bullet$ $Der\ A=
\left\{\left[\begin{array}{lll}x&0&0\\0&y&z\\0&0&0
\end{array}\right]\ |\ x,y,z\in \mathbb{R}\right\}$,

$\bullet$ $Aut\ A= \left\{
\left[\begin{array}{lll}x&0&0\\0&y&z\\0&0&1\end{array}\right]\ |\
x,y, z\in \mathbb{R},\ xy\neq 0 \right\}$,

$\bullet$ the partition $\mathcal{P}_{A}$ of $\mathbb{R}^{3}$,
defined by the lattice of subalgebras of $A$, consists of:

\vspace{1mm} \hspace{5mm} $\diamond$  the singletons covering the
plane $x^{1}Ox^{2}$,

\hspace{5mm} $\diamond$  the half-axes delimited by
 $O$ on the lines $\mathbb{R}u$ for
$u\in \mathcal{I}(A)$ (they cover plane $x^{2}Ox^{3}$ less axis
$Ox^{2}$),

\hspace{5mm} $\diamond$  the quarters  of plane $x^{2}Ox^{3}$,
delimited by axes $Ox^{2}$ and $Ox^{3}$,

\hspace{5mm} $\diamond$ the quarters  of planes passing through
$Ox^{1}$ (less the plane $x^{1}Ox^{2}$), delimited by plane
$x^{2}Ox^{3}$ and axis $Ox^{1}$,

\vspace{1mm}$\bullet$ the partition $\mathcal{P}_{A}$ of $A$
induces a partition on the set of integral curves of the
associated homogeneous quadratic differential system (HQDS)
consisting of:

\vspace{1mm} \hspace{5mm} $\diamond$  the singletons consisting of
singular solutions that cover the plane $x^{1}Ox^{2}$,

\hspace{5mm} $\diamond$  the families of ray solutions lying on
semi-axes of $Ox^{3}$ delimited by $O$,

\hspace{5mm} $\diamond$  the integral curves contained in each
quarter of plane of $x^{2}Ox^{3}$, delimited by axes $Ox^{2}$ and
$Ox^{3}$,

\hspace{5mm} $\diamond$ the integral curves contained in each
quarter of planes passing through $Ox^{1}$ (less the plane
$x^{1}Ox^{2}$), delimited by plane $x^{2}Ox^{3}$ and axis
$Ox^{1}$.

\vspace{3mm}Note that each nonsingular integral curve is
torsion-free.

\vspace{5mm}\textbf{10)}  \emph{Properties of algebras $A$ of type
$A_{2}(\frac{1}{2},\frac{1}{2})$}

\vspace{2mm}$\bullet$ $Ann\ A=\{0\},\ \mathcal{N}(A)=
Span_{\mathbb{R}}\{e_{1}, e_{2}\},\
\mathcal{I}(A)=\{xe_{1}+ye_{2}+e_{3}\ | \ x,y\in\mathbb{R}\}$,

$\bullet$ 1-dimensional subalgebras: $\mathbb{R}u$ for $u\in
\mathcal{N}(A)\cup \mathcal{I}(A)$,

$\bullet$ 2-dimensional subalgebras:
$Span_{\mathbb{R}}\{pe_{1}+qe_{2}+re_{3},\ ae_{1}+be_{2}+ce_{3}\}$
when $rank \left[\begin{array}{lll}p&q&r\\a&b&c
\end{array}\right]=2$ (i.e. each 2-dimensional subspace of $A$ is a subalgebra),

$\bullet$ ideals: $\mathbb{R}e_{1},\ \mathbb{R}e_{2},\
Span_{\mathbb{R}}\{e_{1}, e_{2}\}$,

$\bullet$ $A^{2}=A$,

$\bullet$ $Der\ A=\left\{
\left[\begin{array}{lll}x&y&z\\u&v&w\\0&0&0\end{array}\right]\ |\
x,y,z,u,v,w\in \mathbb{R}\right\}\cong aff\ (2,\mathbb{R})$

$\bullet$ $Aut\ A= \left\{
\left[\begin{array}{lll}x&y&z\\u&v&w\\0&0&1\end{array}\right]\ |\
x,y,z,u,v,w\in \mathbb{R}^{\ast}\right\}\cong Aff\
(2,\mathbb{R})$,

$\bullet$ the partition $\mathcal{P}_{A}$ of $\mathbb{R}^{3}$,
defined by the lattice of subalgebras of $A$, consists of:

\vspace{1mm} \hspace{5mm} $\diamond$  the singletons covering the
plane $x^{1}Ox^{2}$,

\hspace{5mm} $\diamond$ the half-lines delimited by
 $O$ on each line $\mathbb{R}u$ for $u\in \mathcal{I}(A)$ (these lines cover the space less the plane $x^{1}Ox^{2}$),

\vspace{1mm}$\bullet$ the partition $\mathcal{P}_{A}$ of $A$
induces a partition on the set of integral curves of the
associated homogeneous quadratic differential system (HQDS)
consisting of:

\vspace{1mm} \hspace{5mm} $\diamond$  the singletons consisting of
singular solutions that cover the plane $x^{1}Ox^{2}$,

\hspace{5mm} $\diamond$  the ray solutions lying on semi-lines
delimited by $O$ on each line $\mathbb{R}u$ for $u\in
\mathcal{I}(A)$ (these lines cover the space less the plane
$x^{1}Ox^{2}$),

\vspace{3mm} Each nonsingular integral curve lies on a line
$\mathbb{R}u$ for $u\in \mathcal{I}(A)$, so that both its
curvature and torsion tensors vanish.

\vspace{3mm} Note that all idempotents have the same spectrum
$(\frac{1}{2},\frac{1}{2},1)$.

\vspace{5mm}\textbf{A11)}  \emph{Properties of algebras $A$ of
type $A_{2}(0,\beta)$ with} $\beta\notin\{0, \frac{1}{2}\}$

\vspace{3mm} $\bullet$ algebras $A_{2}(0,\beta)$\
($\beta\notin\{0, \frac{1}{2}\}$) and $A_{2}(0,\beta')$\
($\beta'\notin\{0, \frac{1}{2}\}$) are isomorphic if and only if
$\beta=\beta'$,

$\bullet$ $Ann\ A=\mathbb{R}e_{1},\ \mathcal{N}(A)=
Span_{\mathbb{R}}\{e_{1}, e_{2}\},\ \mathcal{I}(A)=\{e_{3}\}$,

$\bullet$ 1-dimensional subalgebras: $\mathbb{R}u$ for $u\in
\mathcal{N}(A)\cup \mathcal{I}(A)$,

$\bullet$ 2-dimensional subalgebras: $Span_{\mathbb{R}}\{e_{1},
e_{2}\},\ Span_{\mathbb{R}}\{e_{1}, e_{3}\}$,
$Span_{\mathbb{R}}\{e_{2}, e_{3}\}$,

$\bullet$ ideals: $\mathbb{R}e_{1},\ \mathbb{R}e_{2},\
Span_{\mathbb{R}}\{e_{1}, e_{2}\},\ Span_{\mathbb{R}}\{e_{2},
e_{3}\}$,

$\bullet$ $A^{2}=Span_{\mathbb{R}}\{e_{2}, e_{3}\}$; $A/A^{2}$ is
the null 1-dimensional algebra; $A$ is a vector direct sum of
ideals,

$\bullet$ $Der\ A=
\left\{\left[\begin{array}{lll}x&0&0\\0&y&0\\0&0&0
\end{array}\right]\ |\ x,y\in \mathbb{R}\right\}$,

$\bullet$ $Aut\ A= \left\{
\left[\begin{array}{lll}x&0&0\\0&y&0\\0&0&1\end{array}\right]\ |\
x,y\in \mathbb{R}^{\ast}\right\}\cong
\mathbb{R}^{\ast}(\cdot)\times\mathbb{R}^{\ast}(\cdot)$.

$\bullet$ the partition $\mathcal{P}_{A}$ of $\mathbb{R}^{3}$,
defined by the lattice of subalgebras of $A$, consists of:

\vspace{1mm} \hspace{5mm} $\diamond$  the singletons covering the
plane $x^{1}Ox^{2}$,

\hspace{5mm} $\diamond$  the half-axes of $Ox^{3}$, delimited by
 $O$,

\hspace{5mm} $\diamond$  the quarters  of plane $x^{1}Ox^{3}$,
delimited by axes $Ox^{1}$ and $Ox^{3}$,

\hspace{5mm} $\diamond$  the quarters  of plane $x^{2}Ox^{3}$,
delimited by axes $Ox^{2}$ and $Ox^{3}$,

\hspace{5mm} $\diamond$ the connected sets in space delimited by
planes $x^{1}Ox^{2},\ x^{1}Ox^{3}$ and $x^{2}Ox^{3}$,

\vspace{1mm}$\bullet$ the partition $\mathcal{P}_{A}$ of $A$
induces a partition on the set of integral curves of the
associated homogeneous quadratic differential system (HQDS)
consisting of:

\vspace{1mm} \hspace{5mm} $\diamond$  the singletons consisting of
singular solutions that cover the plane $x^{1}Ox^{2}$,

\hspace{5mm} $\diamond$  the families of ray solutions lying on
semi-axes of $Ox^{3}$, delimited by $O$,

\hspace{5mm} $\diamond$  the integral curves contained in each
quarters of plane $x^{1}Ox^{3}$, delimited by axes $Ox^{1}$ and
$Ox^{3}$,

\hspace{5mm} $\diamond$  the integral curves contained in each
quarters of plane $x^{2}Ox^{3}$, delimited by axes $Ox^{2}$ and
$Ox^{3}$,

\hspace{5mm} $\diamond$ the integral curves contained in each
 connected sets in space delimited by
planes $x^{1}Ox^{2},\ x^{1}Ox^{3}$ and $x^{2}Ox^{3}$.

\vspace{3mm}\noindent\textbf{Note.} In order to bring together the
classes $A_{2}(\alpha,0)$ and $A_{2}(0,\beta)$ we have ignored, in
case \textbf{11)}, the convention to consider increasing
parameters in algebras of type $A_{2}(a,b)$.

\vspace{5mm}\textbf{A12)}  \emph{Properties of algebras $A$ of
type $A_{2}(\frac{1}{2},\beta)$ with $\beta\notin \{0,
\frac{1}{2}\}$}

\vspace{3mm} $\bullet$ algebras $A_{2}(\frac{1}{2},\beta)\ (
\beta\notin \{0, \frac{1}{2}\})$ and $A_{2}(\frac{1}{2},\beta')\
(\beta'\notin \{0, \frac{1}{2}\})$ are isomorphic if and only if
$\beta=\beta'$,

$\bullet$ $Ann\ A=\{0\},\ \mathcal{N}(A)=
Span_{\mathbb{R}}\{e_{1}, e_{2}\},\ \mathcal{I}(A)=\{xe_{1}+e_{3}\
| \ x\in\mathbb{R}\}$,

$\bullet$ 1-dimensional subalgebras: $\mathbb{R}u$ for $u\in
\mathcal{N}(A)\cup \mathcal{I}(A)$,

$\bullet$ 2-dimensional subalgebras: $Span_{\mathbb{R}}\{e_{1},
e_{3}\}$, $Span_{\mathbb{R}}\{e_{2}, ae_{1}+ce_{3}\}\
(a^{2}+c^{2}\neq 0)$,

$\bullet$ ideals: $\mathbb{R}e_{1},\ \mathbb{R}e_{2},\
Span_{\mathbb{R}}\{e_{1}, e_{2}\}$,

$\bullet$ $A^{2}=A$,

$\bullet$ $Der\ A=\left\{
\left[\begin{array}{lll}x&0&0\\0&y&0\\0&0&0\end{array}\right]\ |\
x,y\in \mathbb{R}\right\}$

$\bullet$ $Aut\ A= \left\{
\left[\begin{array}{lll}x&0&0\\0&y&0\\0&0&1\end{array}\right]\ |\
x,y\in \mathbb{R}^{\ast}\right\}\cong
\mathbb{R}^{\ast}\times\mathbb{R}^{\ast}$,

$\bullet$ the partition $\mathcal{P}_{A}$ of $\mathbb{R}^{3}$,
defined by the lattice of subalgebras of $A$, consists of:

\vspace{1mm} \hspace{5mm} $\diamond$  the singletons covering the
plane $x^{1}Ox^{2}$,

\hspace{5mm} $\diamond$ the half-axes delimited by
 $O$ on each $\mathbb{R}u$ for $u\in \mathcal{I}(A)$ (these lines cover $x^{1}Ox^{3}$ less axis $Ox^{1}$),

\hspace{5mm} $\diamond$ the quarters of space delimited by planes
$x^{1}Ox^{3}$ and $x^{2}Ox^{3}$,

\vspace{1mm}$\bullet$ the partition $\mathcal{P}_{A}$ of $A$
induces a partition on the set of integral curves of the
associated homogeneous quadratic differential system (HQDS)
consisting of:

\vspace{1mm} \hspace{5mm} $\diamond$  the singletons consisting of
singular solutions that cover the plane $x^{1}Ox^{2}$,

\hspace{5mm} $\diamond$  the families of ray solutions lying on
semi-axes delimited by $O$ on each $\mathbb{R}u$ for $u\in
\mathcal{I}(A)$

\hspace{5mm} $\diamond$ the integral curves contained in each
quarter of space delimited by planes $x^{1}Ox^{2}$ and
$x^{1}Ox^{3}$.

\vspace{2mm} All idempotents have the same spectrum:
$\{\beta,\frac{1}{2},1\}$.

\vspace{5mm}\textbf{A13)}  \emph{Properties of algebras $A$ of
type $A_{2}(\alpha,\alpha)$ with $\alpha\notin \{0,
\frac{1}{2}\}$}

\vspace{3mm} $\bullet$ algebras $A_{2}(\alpha,\alpha)\
(\alpha\notin \{0, \frac{1}{2}\})$ and $A_{2}(\alpha',\alpha')\
(\alpha'\notin \{0, \frac{1}{2}\})$ are isomorphic if and only if
$\alpha=\alpha'$,

$\bullet$ $Ann\ A=\{0\},\ \mathcal{N}(A)=
Span_{\mathbb{R}}\{e_{1}, e_{2}\}$,\ $\mathcal{I}(A)=\{e_{3}\}$,

$\bullet$ 1-dimensional subalgebras: $\mathbb{R}u$ for $u\in
\mathcal{N}(A)\cup \mathcal{I}(A)$,

$\bullet$ 2-dimensional subalgebras: $Span_{\mathbb{R}}\{e_{1},
e_{2}\}$, $Span_{\mathbb{R}}\{e_{3}, pe_{1}+qe_{2}\}\
(p^{2}+q^{2}\neq 0)$,

$\bullet$ ideals: $\mathbb{R}e_{1},\ \mathbb{R}e_{2}$,
$Span_{\mathbb{R}}\{e_{1}, e_{2}\}$,

$\bullet$ $A^{2}=A$,

$\bullet$ $Der\ A=\left\{
\left[\begin{array}{lll}x&y&0\\z&v&0\\0&0&0\end{array}\right]\ |\
x,y,z,v\in \mathbb{R}\right\}\cong g\ell(2,\mathbb{R})$,

$\bullet$ $Aut\ A= \left\{
\left[\begin{array}{lll}x&y&0\\z&v&0\\0&0&1\end{array}\right]\ |\
x,y,z,v\in \mathbb{R},\ xv-yz\neq 0\right\}\cong
GL(2,\mathbb{R})$,

$\bullet$ the partition $\mathcal{P}_{A}$ of $\mathbb{R}^{3}$,
defined by the lattice of subalgebras of $A$, consists of:

\vspace{1mm} \hspace{5mm} $\diamond$  the singletons covering the
plane $x^{1}Ox^{2}$,

\hspace{5mm} $\diamond$ the half-lines delimited by
 $O$ on $Ox^{3}$,

\hspace{5mm} $\diamond$ the quarters delimited by axis $Ox^{3}$
and plane $x^{1}Ox^{2}$ on each plane passing through  $Ox^{3}$,

\vspace{1mm}$\bullet$ the partition $\mathcal{P}_{A}$ of $A$
induces a partition on the set of integral curves of the
associated homogeneous quadratic differential system (HQDS)
consisting of:

\vspace{1mm} \hspace{5mm} $\diamond$  the singletons consisting of
singular solutions that cover the plane $x^{1}Ox^{2}$,

\hspace{5mm} $\diamond$  the families of ray solutions lying on
semi-axes delimited by $O$ on axis $Ox^{3}$,

\hspace{5mm} $\diamond$ the integral curves contained in the
quarters delimited by axis $Ox^{3}$ and plane $x^{1}Ox^{2}$ on
each plane passing through $Ox^{3}$.

\vspace{5mm}\textbf{A14)}  \emph{Properties of algebras $A$ of
type $A_{2}(\alpha,\beta)$ with $\alpha,\ \beta\notin \{0,
\frac{1}{2}\}\ and\ \alpha<\beta$}

\vspace{3mm} $\bullet$ algebras $A_{2}(\alpha,\beta)\
(\alpha<\beta\ and \ \alpha,\ \beta\notin \{0, \frac{1}{2}\})$ and
$A_{2}(\alpha',\beta')\ (\alpha'<\beta'\ and \ \alpha',\
\beta'\notin \{0, \frac{1}{2}\})$ are isomorphic if and only if
$\alpha=\alpha'\ \ and\ \ \beta=\beta'$,

$\bullet$ $Ann\ A=\{0\},\ \mathcal{N}(A)=
Span_{\mathbb{R}}\{e_{1}, e_{2}\}$,\ $\mathcal{I}(A)=\{e_{3}\}$,

$\bullet$ 1-dimensional subalgebras: $\mathbb{R}u$ for $u\in
\mathcal{N}(A)\cup \mathcal{I}(A)$,

$\bullet$ 2-dimensional subalgebras: $Span_{\mathbb{R}}\{e_{1},
e_{2}\}$, $Span_{\mathbb{R}}\{e_{1}, e_{3}\}$,
$Span_{\mathbb{R}}\{e_{2}, e_{3}\}$,

$\bullet$ ideals: $\mathbb{R}e_{1},\ \mathbb{R}e_{2}$,
$Span_{\mathbb{R}}\{e_{1}, e_{2}\}$,

$\bullet$ $A^{2}=A$,

$\bullet$ $Der\ A=\left\{
\left[\begin{array}{lll}x&0&0\\0&y&0\\0&0&0\end{array}\right]\ |\
x,y\in \mathbb{R}\right\}$

$\bullet$ $Aut\ A= \left\{
\left[\begin{array}{lll}x&0&0\\0&y&0\\0&0&1\end{array}\right]\ |\
x,y\in \mathbb{R}^{\ast}\right\}\cong \mathbb{R}^{\ast}\times
\mathbb{R}^{\ast}$,

$\bullet$ the partition $\mathcal{P}_{A}$ of $\mathbb{R}^{3}$,
defined by the lattice of subalgebras of $A$, consists of:

\vspace{1mm} \hspace{5mm} $\diamond$  the singletons covering the
plane $x^{1}Ox^{2}$,

\hspace{5mm} $\diamond$ the half-axes delimited by
 $O$ on $Ox^{3}$,

\hspace{5mm} $\diamond$ the quarters  of plane $x^{1}Ox^{3}$,
delimited by axes $Ox^{1}$ and $Ox^{3}$,

\hspace{5mm} $\diamond$  the quarters  of plane $x^{2}Ox^{3}$,
delimited by axes $Ox^{2}$ and $Ox^{3}$,

\hspace{5mm} $\diamond$ the connected components of space
delimited by planes $x^{1}Ox^{2}$, $x^{1}Ox^{3}$ and
$x^{2}Ox^{3}$,

\vspace{1mm}$\bullet$ the partition $\mathcal{P}_{A}$ of $A$
induces a partition of the set of integral curves of the
associated homogeneous quadratic differential system (HQDS)
consisting of:

\vspace{1mm} \hspace{5mm} $\diamond$  the singletons consisting of
singular solutions that cover the axes $Ox^{1}$ and $Ox^{2}$,

\hspace{5mm} $\diamond$  the ray solutions lying on semi-axes of
$Ox^{3}$, delimited by $O$,

\hspace{5mm} $\diamond$  the integral curves contained in each
quarter of the plane of $x^{1}Ox^{3}$, delimited by axes $Ox^{1}$
and $Ox^{3}$,

\hspace{5mm} $\diamond$  the integral curves contained in each
quarter of the plane of $x^{2}Ox^{3}$, delimited by axes $Ox^{2}$
and $Ox^{3}$,

\hspace{5mm} $\diamond$ the integral curves contained in each
connected component of space delimited by planes $x^{1}Ox^{2}$,
$x^{1}Ox^{3}$ and $x^{2}Ox^{3}$.

\vspace{3mm} \emph{Subcase (ii)\ j=0,\ $ n\neq 0$}

\vspace{3mm} There exists a basis such that the multiplication
table of algebra becomes

$$\begin{array}{llll}
  \textbf{Table T3}\hspace{6mm} & \hspace{5mm} e_{1}^{2}=0  &\hspace{5mm} e_{2}^{2}=0&\hspace{5mm} e_{3}^{2}=0 \\
   &\hspace{5mm}  e_{1}e_{2}=e_{3}&\hspace{5mm}  e_{1}e_{3}=\alpha e_{1}  &\hspace{5mm}  e_{2}e_{3}=\beta e_{2}
\end{array}$$
with $\alpha,\ \beta\in \mathbb{R}$. Let us denote by
$A_{3}(\alpha, \beta)$ any algebra of type $\textbf{T3}$.

\begin{proposition}
The algebras $A_{3}(\alpha, \beta)$ and $A_{3}( \beta, \alpha)$
are isomorphic.
\end{proposition}
Consequently, in the following we restrict our interest to the
case $\alpha\leq \beta$.

\begin{proposition} Every algebra $A$ of type $A_{3}(\alpha,\beta)$ has:
$$Ann\ A=\left\{\begin{array}{lll}\mathbb{R}e_{3}&if&\alpha=\beta=0\\
\{0\}&if&\alpha^{2}+\beta^{2}\neq 0 \end{array} \right.$$
$$\mathcal{N}(A)=\left\{\begin{array}{lll}Span_{\mathbb{R}}\{e_{1}, e_{3}\}\cup Span_{\mathbb{R}}\{e_{2},
e_{3}\}&if&\alpha=\beta=0\\

\mathbb{R}e_{1}\cup Span_{\mathbb{R}}\{e_{2},
e_{3}\}&if&\alpha\neq 0,\ \beta=0\\

\mathbb{R}e_{2}\cup Span_{\mathbb{R}}\{e_{1},
e_{3}\}&if&\alpha = 0,\ \beta\neq 0\\

\mathbb{R}e_{1}\cup \mathbb{R}e_{2}\cup
\mathbb{R}e_{3}&if&\alpha\beta\neq 0
\end{array}\right.$$

$$\mathcal{I}(A)=\left\{\begin{array}{clc}\emptyset&if&\alpha=0\bigvee \beta=0\bigvee \alpha\neq \beta\\
\{xe_{1}+\frac{1}{4\alpha x}e_{2}+\frac{1}{2\alpha}e_{3}\ |\ x\in
\mathbb{R}^{\ast}\}&if&\alpha=\beta.

\end{array}\right.$$
\end{proposition}

\vspace{5mm}\textbf{A15)} \emph{Properties of algebras $A$ of type
$A_{3}(0,0)$}

\vspace{3mm}It is suitable to use the change of bases
$(e_{1},e_{2},e_{3})\rightarrow (e_{3},e_{2},e_{1})$.

\vspace{2mm}$\bullet$ $Ann\ A=\mathbb{R}e_{1},\ \mathcal{N}(A)=
Span_{\mathbb{R}}\{e_{1}, e_{2}\}\cup Span_{\mathbb{R}}\{e_{1},
e_{3}\},\ \mathcal{I}(A)=\emptyset$,

$\bullet$ 1-dimensional subalgebras: $\mathbb{R}u$ for $u\in
\mathcal{N}(A)$,

$\bullet$ 2-dimensional subalgebras: $Span_{\mathbb{R}}\{e_{1},
be_{2}+ce_{3}\}\ (b^{2}+c^{2}\neq 0)$,

$\bullet$ ideals: $\mathbb{R}e_{1},\ Span_{\mathbb{R}}\{e_{1},
be_{2}+ce_{3}\}\ (b^{2}+c^{2}\neq 0)$,

$\bullet$ $A^{2}=\mathbb{R}e_{1}$; $A/A^{2}$ is a 2-dimensional
null algebra,

$\bullet$ $Der\ A=
\left\{\left[\begin{array}{lll}x+y&z&u\\0&y&0\\0&0&x
\end{array}\right]\ |\ x,y,z,u\in \mathbb{R}\right\}$,

$\bullet$ $Aut\ A=H\cup J\cdot H$, where

$$H= \left\{
\left[\begin{array}{lll}xy&z&v\\0&y&0\\0&0&x\end{array}\right]\ |\
x,y,z,v\in \mathbb{R},\ xy\neq 0\right\}\ \ \ \
J=\left[\begin{array}{lll}1&0&0\\0&0&1\\0&1&0\end{array}\right];$$
$H$ is a normal divisor of $Aut\ A$,

$\bullet$ the partition $\mathcal{P}_{A}$ of $\mathbb{R}^{3}$,
defined by the lattice of subalgebras of $A$, consists of:

\vspace{1mm} \hspace{5mm} $\diamond$  the singletons covering the
planes $x^{1}Ox^{2}$ and  $x^{1}Ox^{3}$,

\hspace{5mm} $\diamond$ the half-planes delimited by axis $Ox^{1}$
on each plane passing through  $Ox^{1}$ less the planes
$x^{1}Ox^{2}$ and $x^{1}Ox^{3}$,

\vspace{1mm}$\bullet$ the partition $\mathcal{P}_{A}$ of $A$
induces a partition on the set of integral curves of the
associated homogeneous quadratic differential system (HQDS)
consisting of:

\vspace{1mm} \hspace{5mm} $\diamond$  the singletons consisting of
singular solutions that cover the planes $x^{1}Ox^{2}$ and
$x^{1}Ox^{3}$,

\hspace{5mm} $\diamond$ the nonsingular integral curves contained
in each half-plane delimited by axis $Ox^{1}$ on each plane
passing through  $Ox^{1}$ less the planes $x^{1}Ox^{2}$ and
$x^{1}Ox^{3}$.

\vspace{2mm}Since each nonsingular integral curve lies on a plane
passing through $Ox^{1}$ it has a zero torsion tensor. In
addition, $A^{2}=\mathbb{R}e_{1}$ assures that each nonsingular
integral curve has the curvature tensors zero, so that each
nonsingular integral curve lies on a line parallel to $Ox^{1}$.

\vspace{3mm}\emph{Case} $\alpha^{2}+\beta^{2}\neq 0$

\vspace{3mm}The class of algebras $A(\alpha,\beta)$ with
$\alpha^{2}+\beta^{2}\neq 0$ decomposes naturally in accordance
with the next conditions:

$$\begin{array}{lll}(i)&\alpha\neq 0,\ \beta=0 \\
(ii)&\alpha =0,\ \beta\neq 0 \\
(iii)&\alpha\beta\neq 0& \alpha=\beta \\
(iv)&\alpha\beta\neq 0& \alpha\neq \beta, \ \alpha < \beta.
\end{array}$$

\begin{proposition} The following assertions hold:

\vspace{2mm} (i)\  \ each algebra $A_{3}(\alpha, 0)$ with
$\alpha\neq 0$ is isomorphic to algebra $A_{3}(1,0)\cong
A_{3}(0,1)$,

(ii)\  each algebra $A_{3}(0,\beta)$ with $\beta\neq 0$ is
isomorphic to algebra $A_{3}(0,1)\cong A_{3}(1,0)$,

(iii) each algebra $A_{3}(\alpha, \beta)$ with $\alpha=\beta$ is
isomorphic to algebra $A_{3}(1,1)$,

(iv)\  each algebra $A_{3}(\alpha,\beta)$ with $\alpha\beta\neq 0$
and $\alpha<\beta$ is isomorphic to algebra $A_{3}(1,\beta)$.
\end{proposition}

\vspace{3mm}\textbf{A16)} \emph{Properties of algebra
$A=A_{3}(0,1)$}

\vspace{3mm}$\bullet$ $Ann\ A=\{0\},\
\mathcal{N}(A)=\mathbb{R}e_{2}\cup Span_{\mathbb{R}}\{e_{1},
e_{3}\},\ \mathcal{I}(A)=\emptyset$,

$\bullet$ 1-dimensional subalgebras: $\mathbb{R}u$ for $u\in
\mathcal{N}(A)$,

$\bullet$ 2-dimensional subalgebras: $Span_{\mathbb{R}}\{e_{1},
e_{3}\},\ Span_{\mathbb{R}}\{e_{2}, e_{3}\}$,

$\bullet$ ideals: $Span_{\mathbb{R}}\{e_{2}, e_{3}\}$,

$\bullet$ $A^{2}=Span_{\mathbb{R}}\{e_{2}, e_{3}\}$; $A/A^{2}$ is
a 1-dimensional null algebra,

$\bullet$ $Der\ A=\mathbb{R}D$,

$\bullet$ $Aut\ A=\left\{
\left[\begin{array}{lll}x&0&0\\0&\frac{1}{x}&0\\0&0&1\end{array}\right]\
|\ x\in \mathbb{R}^{\ast}\right\}\cong\mathbb{R}^{\ast}(\cdot)$

$\bullet$ the partition $\mathcal{P}_{A}$ of $\mathbb{R}^{3}$,
defined by the lattice of subalgebras of $A$, consists of:

\vspace{1mm} \hspace{5mm} $\diamond$  the singletons covering the
axis $Ox^{2}$ and the plane $x^{1}Ox^{3}$,

\hspace{5mm} $\diamond$  the connected sets delimited by axis
$Ox^{2}$ and the plane $x^{1}Ox^{3}$,

\vspace{1mm}$\bullet$ the partition $\mathcal{P}_{A}$ of $A$
induces a partition on the set of integral curves of the
associated homogeneous quadratic differential system (HQDS)
consisting of:

\vspace{1mm} \hspace{5mm} $\diamond$  the singletons consisting of
singular solutions that cover the axis $Ox^{2}$ and the plane
$x^{1}Ox^{3}$,

\hspace{5mm} $\diamond$  the integral curves contained in each
connected sets delimited by axis $Ox^{2}$ and the planes
$x^{1}Ox^{3}$.

\vspace{3mm} Since algebra $A$ is solvable, the nonsingular
integral curves are torsion-free.

\vspace{3mm}\textbf{A17)} \emph{Properties of algebra
$A=A_{3}(1,1)$}

\vspace{3mm}$\bullet$ $Ann\ A=\{0\},\
\mathcal{N}(A)=\mathbb{R}e_{1}\cup \mathbb{R}e_{2}\cup
\mathbb{R}e_{3},\
\mathcal{I}(A)=\{xe_{1}+\frac{1}{4x}e_{2}+\frac{1}{2}e_{3}\ |\
x\in \mathbb{R}^{\ast}\}$,

$\bullet$ 1-dimensional subalgebras: $\mathbb{R}u$ for $u\in
\mathcal{N}(A)\cup \mathcal{I}(A)$,

$\bullet$ 2-dimensional subalgebras: $Span_{\mathbb{R}}\{e_{3},
ae_{1}+be_{2}\}$,

$\bullet$ ideals: none,

$\bullet$ $A^{2}=A$,

$\bullet$ $Der\ A=\mathbb{R}D$,

$\bullet$ $Aut\ A=H\cup JH$ where

$$H=\left\{
\left[\begin{array}{lll}x&0&0\\0&\frac{1}{x}&0\\0&0&1\end{array}\right]\
|\ x\in \mathbb{R}^{\ast}\right\},\ \
J=\left[\begin{array}{lll}0&1&0\\1&0&0\\0&0&1\end{array}\right],$$
i.e. $Aut\ A\cong \mathbb{R}^{\ast}\times \{-1,1\}$,

$\bullet$ the partition $\mathcal{P}_{A}$ of $\mathbb{R}^{3}$,
defined by the lattice of subalgebras of $A$, consists of:

\vspace{1mm} \hspace{5mm} $\diamond$  the singletons covering the
axes $Ox^{1}$,  $Ox^{2}$ and $Ox^{3}$,

\vspace{1mm} \hspace{5mm} $\diamond$  the half-lines delimited by
$O$ on each line passing through $O$ and directed by an
idempotent,

\hspace{5mm} $\diamond$  the connected sets delimited by axis
$Ox^{3}$ and the line $\mathbb{R}E$ directed by an idempotent
$E=\omega e_{1}+\frac{1}{4\omega}e_{2}+\frac{1}{2}e_{3}$ on the
plane parallel to $E$ passing through $Ox^{3}$,

\hspace{5mm} $\diamond$  the half-planes delimited by axis
$Ox^{3}$ on any plane  passing through $Ox^{3}$ which is not
parallel to any idempotent (i.e. with $x^{1}x^{2}<0$),

\vspace{1mm}$\bullet$ the partition $\mathcal{P}_{A}$ of $A$
induces a partition on the set of integral curves of the
associated homogeneous quadratic differential system (HQDS)
consisting of:

\vspace{1mm} \hspace{5mm} $\diamond$  the singletons consisting of
singular solutions that cover the axes $Ox^{1}$,  $Ox^{2}$ and
$Ox^{3}$,

\hspace{5mm} $\diamond$  the families of ray solutions contained
in each half-line delimited by $O$ on each line passing through
$O$ and directed by an idempotent,

\hspace{5mm} $\diamond$  the integral curves contained in the
connected sets delimited by axis $Ox^{3}$ and the line
$\mathbb{R}E$ directed by an idempotent $E=\omega
e_{1}+\frac{1}{4\omega}e_{2}+\frac{1}{2}e_{3}$ on the plane
parallel to $E$ passing through $Ox^{3}$,

\hspace{5mm} $\diamond$  the integral curves contained in the
half-planes delimited by axis $Ox^{3}$ on any plane  passing
through $Ox^{3}$ which is not parallel to any idempotent (i.e.
with $x^{1}x^{2}<0$).

\vspace{3mm} All nonsingular integral curves are torsion-free.

\vspace{3mm}\textbf{A18)} \emph{Properties of algebra $A$ of type
$A_{3}(1, \beta)$ with $\beta>1$}

\vspace{3mm}$\bullet$ the algebras $A_{3}(1, \beta)$ $(\beta>1)$
and $A_{3}(1, \beta')$ $(\beta'>1)$ are isomorphic if and only if
$\beta=\beta'$,

$\bullet$ $Ann\ A=\{0\},\ \mathcal{N}(A)=\mathbb{R}e_{1}\cup
\mathbb{R}e_{2}\cup \mathbb{R}e_{3},\ \mathcal{I}(A)=\emptyset$,

$\bullet$ 1-dimensional subalgebras: $\mathbb{R}e_{1},\
\mathbb{R}e_{2},\ \mathbb{R}e_{3}$,

$\bullet$ 2-dimensional subalgebras: $Span_{\mathbb{R}}\{e_{1},
e_{3}\},\ Span_{\mathbb{R}}\{e_{2}, e_{3}\}$,

$\bullet$ ideals: none,

$\bullet$ $A^{2}=A$,

$\bullet$ $Der\ A=\mathbb{R}D$,

$\bullet$ $Aut\ A=H\cup JH$ where

$$H=\left\{
\left[\begin{array}{lll}x&0&0\\0&\frac{1}{x}&0\\0&0&1\end{array}\right]\
|\ x\in \mathbb{R}^{\ast}\right\},\ \
J=\left[\begin{array}{lll}0&1&0\\1&0&0\\0&0&1\end{array}\right],$$
i.e. $Aut\ A\cong \mathbb{R}^{\ast}\times \{-1,1\}$,

$\bullet$ the partition $\mathcal{P}_{A}$ of $\mathbb{R}^{3}$,
defined by the lattice of subalgebras of $A$, consists of:

\vspace{1mm} \hspace{5mm} $\diamond$  the singletons covering the
axes $Ox^{1}$,  $Ox^{2}$ and $Ox^{3}$,

\vspace{1mm} \hspace{5mm} $\diamond$ the quarters of plane
$x^{1}Ox^{3}$ delimited by axes $Ox^{1}$ and $Ox^{3}$,

\vspace{1mm} \hspace{5mm} $\diamond$ the quarters of plane
$x^{2}Ox^{3}$ delimited by axes $Ox^{2}$ and $Ox^{3}$,

\vspace{1mm} \hspace{5mm} $\diamond$ the quarters of space
delimited by planes $x^{1}Ox^{3}$ and $x^{2}Ox^{3}$,

\vspace{1mm}$\bullet$ the partition $\mathcal{P}_{A}$ of $A$
induces a partition on the set of integral curves of the
associated homogeneous quadratic differential system (HQDS)
consisting of:

\vspace{1mm} \hspace{5mm} $\diamond$  the singletons consisting of
singular solutions that cover the axes $Ox^{1}$,  $Ox^{2}$ and
$Ox^{3}$,

\vspace{1mm} \hspace{5mm} $\diamond$  the integral curves
contained in the quarters of plane $x^{1}Ox^{3}$ delimited by axes
$Ox^{1}$ and $Ox^{3}$,

\vspace{1mm} \hspace{5mm} $\diamond$  the integral curves
contained in the quarters of plane $x^{2}Ox^{3}$ delimited by axes
$Ox^{2}$ and $Ox^{3}$,

\vspace{1mm} \hspace{5mm} $\diamond$  the integral curves
contained in the quarters of space delimited by planes
$x^{1}Ox^{3}$ and $x^{2}Ox^{3}$.

\vspace{5mm} \emph{Subcase (iii)\ \ $j=n=0$}

\vspace{3mm} There exists a basis such that the multiplication
table of algebra becomes

$$\begin{array}{llll}
  \textbf{Table T4}\hspace{6mm} & \hspace{5mm} e_{1}^{2}=0  &\hspace{5mm} e_{2}^{2}=0&\hspace{5mm} e_{3}^{2}=0 \\
   &\hspace{5mm}  e_{1}e_{2}=0&\hspace{5mm}  e_{1}e_{3}=\alpha e_{1}  &\hspace{5mm}  e_{2}e_{3}=\beta e_{2}
\end{array}$$
with $\alpha,\ \beta\in \mathbb{R}$. Let us denote by
$A_{4}(\alpha,\beta)$ each algebra defined by multiplication table
\textbf{T4}.

\begin{proposition} Algebras $A_{4}(\alpha, \beta)$ and $A_{4}(\beta, \alpha)$ are isomorphic.
\end{proposition}
Consequently, in the following we will be interested only in the
case $\alpha \leq \beta$.

\vspace{3mm} We have to consider the next mutually exclusive
situations:

\vspace{1mm} (i)  \ \ \  $\alpha\neq 0,\ \ \beta= 0$,

(ii)  \ \ $\alpha =0,\ \ \beta\neq 0$,

(iii)   \ $\alpha\beta\neq 0,\ \ \alpha=\beta$,

(iv) \ \ $\alpha\beta\neq 0$ and $\alpha\neq\beta$.
\begin{proposition} The following assertions hold:

\vspace{2mm}(i) each algebra $A_{4}(\alpha, 0)$ with $\alpha\neq
0$ is isomorphic to algebra $A_{4}(1,0)$,

(ii) each algebra $A_{4}(0,\beta)$ with $\beta\neq 0$ is
isomorphic to algebra $A_{4}(0,1)\cong A_{4}(1,0)$,

(iii) each algebra $A_{4}(\alpha, \beta)$ with $\alpha=\beta$ is
isomorphic to algebra $A_{4}(1,1)$,

(ii) each algebra $A_{4}(\alpha,\beta)$ with $\alpha\beta\neq 0$
and $\alpha\neq\beta$
 is isomorphic to algebra $A_{4}(1,\beta)$.
\end{proposition}

\vspace{3mm} \textbf{A19)} \emph{Properties of algebra
$A=A_{4}(0,1)$}

\vspace{3mm}$\bullet$ $Ann\ A=\mathbb{R}e_{1},\
\mathcal{N}(A)=Span_{\mathbb{R}}\{e_{1}, e_{2}\}\cup
Span_{\mathbb{R}}\{e_{1}, e_{3}\},\ \mathcal{I}(A)=\emptyset$,

$\bullet$ 1-dimensional subalgebras: $\mathbb{R}u$ for $u\in
\mathcal{N}(A)$,

$\bullet$ 2-dimensional subalgebras: $Span_{\mathbb{R}}\{e_{1},
e_{3}\},\ Span_{\mathbb{R}}\{e_{2}, ae_{1}+be_{3}\}\
(a^{2}+b^{2}\neq 0)$,

$\bullet$ ideals:  $\mathbb{R}e_{1},\ \mathbb{R}e_{2},\
Span_{\mathbb{R}}\{e_{2}, ae_{1}+be_{3}\} \ (a^{2}+b^{2}\neq 0)$,

$\bullet$ $A^{2}=\mathbb{R}e_{2}$; $A/A^{2}$ is a 2-dimensional
null algebra,; $A/A^{2}$ is a 2-dimensional null algebra,

$\bullet$ $Der\ A=\left\{
\left[\begin{array}{lll}x&0&y\\0&z&0\\0&0&0\end{array}\right]\ |\
x,y,z\in \mathbb{R}\right\}$,

$\bullet$ $Aut\ A=\left\{
\left[\begin{array}{lll}x&0&y\\0&z&0\\0&0&1\end{array}\right]\ |\
x,y,z\in \mathbb{R},\ xz\neq 0\right\},$

$\bullet$ the partition $\mathcal{P}_{A}$ of $\mathbb{R}^{3}$,
defined by the lattice of subalgebras of $A$, consists of:

\vspace{1mm} \hspace{5mm} $\diamond$  the singletons covering the
planes $x^{1}Ox^{2}$ and $x^{1}Ox^{3}$,

\vspace{1mm} \hspace{5mm} $\diamond$ the quarters of space
delimited by planes $x^{1}Ox^{2}$ and $x^{1}Ox^{3}$,

\vspace{1mm}$\bullet$ the partition $\mathcal{P}_{A}$ of $A$
induces a partition on the set of integral curves of the
associated homogeneous quadratic differential system (HQDS)
consisting of:

\vspace{1mm} \hspace{5mm} $\diamond$  the singletons consisting of
singular solutions that cover the planes $x^{1}Ox^{2}$ and
$x^{1}Ox^{3}$,

\vspace{1mm} \hspace{5mm} $\diamond$  the integral curves
contained in the quarters of space delimited by planes
$x^{1}Ox^{2}$ and $x^{1}Ox^{3}$.

\vspace{3mm}\textbf{A20)} \emph{Properties of algebra
$A=A_{4}(1,1)$}

\vspace{3mm}$\bullet$ $Ann\ A=\{0\},\
\mathcal{N}(A)=Span_{\mathbb{R}}\{e_{1}, e_{2}\}\cup
\mathbb{R}e_{3},\ \mathcal{I}(A)=\emptyset$,

$\bullet$ 1-dimensional subalgebras: $\mathbb{R}u$ for $u\in
\mathcal{N}(A)$,

$\bullet$ 2-dimensional subalgebras: $Span_{\mathbb{R}}\{e_{1},
e_{2}\},\ Span_{\mathbb{R}}\{e_{3}, ae_{1}+be_{2}\}\
(a^{2}+b^{2}\neq 0)$,

$\bullet$ ideals:  $\mathbb{R}e_{1},\ \mathbb{R}e_{2},\
Span_{\mathbb{R}}\{e_{1}, e_{2}\}$,

$\bullet$ $A^{2}=Span_{\mathbb{R}}\{e_{1}, e_{2}\}$; $A/A^{2}$ is
a 1-dimensional null algebra,,

$\bullet$ $Der\ A=\left\{
\left[\begin{array}{lll}x&y&0\\z&u&0\\0&0&0\end{array}\right]\ |\
x,y,z,u\in \mathbb{R}\right\}\cong g\ell(2,\mathbb{R})$,

$\bullet$ $Aut\ A=\left\{
\left[\begin{array}{lll}x&y&0\\z&u&0\\0&0&1\end{array}\right]\ |\
x,y,z,u\in \mathbb{R},\ xu-yz\neq 0\right\}$, i.e. $Aut\ A\cong
GL(2,\mathbb{R})$,

$\bullet$ the partition $\mathcal{P}_{A}$ of $\mathbb{R}^{3}$,
defined by the lattice of subalgebras of $A$, consists of:

\vspace{1mm} \hspace{5mm} $\diamond$  the singletons covering the
plane $x^{1}Ox^{2}$ and axis $Ox^{3}$,

\vspace{1mm} \hspace{5mm} $\diamond$ the quarters of plane
delimited by axis $Ox^{3}$ and  $x^{1}Ox^{2}$ on each plane
passing through $Ox^{3}$,

\vspace{1mm}$\bullet$ the partition $\mathcal{P}_{A}$ of $A$
induces a partition on the set of integral curves of the
associated homogeneous quadratic differential system (HQDS)
consisting of:

\vspace{1mm} \hspace{5mm} $\diamond$  the singletons consisting of
singular solutions that cover the plane $x^{1}Ox^{2}$ and the axis
$Ox^{3}$,

\vspace{1mm} \hspace{5mm} $\diamond$  the integral curves
contained in the quarters of plane delimited by axis $Ox^{3}$ and
$x^{1}Ox^{2}$ on each plane passing through $Ox^{3}$.

\vspace{3mm} Consequently, each nonsingular integral curve is
torsion-free.

\vspace{3mm}\textbf{A21)} \emph{Properties of algebra $A$ of type
$A_{4}(1,\beta)\ (\beta\notin\{0,1\})$}

\vspace{3mm}$\bullet$ the algebras $A_{4}(1, \beta)$
$(\beta\notin\{0,1\})$ and $A_{4}(1, \beta')$
$(\beta\notin\{0,1\})$ are isomorphic if and only if
$\beta=\beta'$,

$\bullet$ $Ann\ A=\{0\},\ \mathcal{N}(A)=Span_{\mathbb{R}}\{e_{1},
e_{2}\}\cup \mathbb{R}e_{3},\ \mathcal{I}(A)=\emptyset$,

$\bullet$ 1-dimensional subalgebras: $\mathbb{R}u$ for $u\in
\mathcal{N}(A)$,

$\bullet$ 2-dimensional subalgebras: $Span_{\mathbb{R}}\{e_{1},
e_{2}\},\ Span_{\mathbb{R}}\{e_{1}, e_{3}\}, \ \
Span_{\mathbb{R}}\{e_{2}, e_{3}\}$,

$\bullet$ ideals:  $\mathbb{R}e_{1},\ \mathbb{R}e_{2},\
Span_{\mathbb{R}}\{e_{1}, e_{2}\}$,

$\bullet$ $A^{2}=Span_{\mathbb{R}}\{e_{1}, e_{2}\}$; $A/A^{2}$ is
a 1-dimensional null algebra,

$\bullet$ $Der\ A=\left\{
\left[\begin{array}{lll}x&0&0\\0&y&0\\0&0&0\end{array}\right]\ |\
x,y\in \mathbb{R}\right\}$,

$\bullet$ $Aut\ A=\left\{
\left[\begin{array}{lll}x&0&0\\0&y&0\\0&0&1\end{array}\right]\ |\
x,y\in \mathbb{R}^{\ast}\right\}\cong
\mathbb{R}^{\ast}\times\mathbb{R}^{\ast}$,

$\bullet$ the partition $\mathcal{P}_{A}$ of $\mathbb{R}^{3}$,
defined by the lattice of subalgebras of $A$, consists of:

\vspace{1mm} \hspace{5mm} $\diamond$  the singletons covering the
plane $x^{1}Ox^{2}$ and axis $Ox^{3}$,

\vspace{1mm} \hspace{5mm} $\diamond$ the quarters delimited by
axes $Ox^{1}$ and $Ox^{3}$ on the plane $x^{1}Ox^{3}$,

\vspace{1mm} \hspace{5mm} $\diamond$ the quarters delimited by
axes $Ox^{2}$ and $Ox^{3}$ on the plane $x^{2}Ox^{3}$,

\vspace{1mm} \hspace{5mm} $\diamond$ the connected sets of space
$\mathbb{R}^{3}$ delimited by coordinate planes $x^{1}Ox^{2}$,
$x^{1}Ox^{3}$ and $x^{2}Ox^{3}$,

\vspace{1mm}$\bullet$ the partition $\mathcal{P}_{A}$ of $A$
induces a partition on the set of integral curves of the
associated homogeneous quadratic differential system (HQDS)
consisting of:

\vspace{1mm} \hspace{5mm} $\diamond$  the singletons consisting of
singular solutions that cover the plane $x^{1}Ox^{2}$ and the axis
$Ox^{3}$,

\vspace{1mm} \hspace{5mm} $\diamond$  the integral curves
contained in the quarters of plane $x^{1}Ox^{3}$ delimited by axes
$Ox^{1}$ and $Ox^{3}$,

\vspace{1mm} \hspace{5mm} $\diamond$  the integral curves
contained in the quarters of plane $x^{2}Ox^{3}$ delimited by axes
$Ox^{2}$ and $Ox^{3}$,

\vspace{1mm} \hspace{5mm} $\diamond$  the integral curves
contained in the connected sets of space $\mathbb{R}^{3}$
delimited by coordinate planes $x^{1}Ox^{2}$, $x^{1}Ox^{3}$ and
$x^{2}Ox^{3}$.

\vspace{5mm} \textbf{2) Case} $Spec\ D=(1,1,0)$

\vspace{3mm}Algebra $A(\cdot)$ has the next multiplication table:
$$\begin{array}{llll}
  \textbf{Table T}\hspace{6mm} & \hspace{5mm} e_{1}^{2}=0  &\hspace{5mm} e_{2}^{2}=0&\hspace{5mm} e_{3}^{2}=e_{3} \\
   &\hspace{5mm}  e_{1}e_{2}= 0 &\hspace{5mm}  e_{1}e_{3}=\alpha e_{1}+\beta e_{2} &\hspace{5mm}  e_{2}e_{3}= \gamma e_{1}+\delta e_{2}
\end{array}$$
with $\alpha,\ \beta,\ \gamma,\ \delta\in \mathbb{R}$. It is
natural to consider the next complementary subcases:

\vspace{1mm}\ I)\ $\alpha\delta-\beta\gamma= 0$,

II) $\alpha\delta-\beta\gamma\neq 0$.

\vspace{3mm} In its turn, Subcase I is naturally divided into two
disjoint parts:

\textbf{I}$_{1}$ $\alpha\delta-\beta\gamma= 0,\ \
\alpha\beta\gamma\delta\neq 0$,

\textbf{I}$_{2}$ $\alpha\delta-\beta\gamma= 0,\ \
\alpha\beta\gamma\delta= 0$.

\vspace{3mm} \textbf{I$_{1}$.}\ \
 $\alpha\delta-\beta\gamma= 0,\ \ \alpha\beta\gamma\delta\neq 0$

\vspace{3mm} The algebra $A(\cdot)$ has $Ann\ A=\mathbb{R}(\gamma
e_{1}-\alpha e_{2})$. We have to consider the next complementary
situations:

$$\textbf{I$_{11}$}.\ \ \alpha^{2}+\beta\gamma\neq 0,\ \ \ I_{12}\ \ \alpha^{2}+\beta\gamma= 0.$$

\vspace{3mm}\textbf{I}$_{11}$ By using basis $(\gamma e_{1}-\alpha
e_{2},\ \alpha e_{1}+\beta e_{2},\ e_{3})$ with
$\alpha^{2}+\beta\gamma\neq 0$ the multiplication table of algebra
becomes:
$$\begin{array}{llll}
  \textbf{Table $TI_{11}$}\hspace{6mm} & \hspace{5mm} e_{1}^{2}=0  &\hspace{5mm} e_{2}^{2}=0&\hspace{5mm} e_{3}^{2}=e_{3} \\
   &\hspace{5mm}  e_{1}e_{2}=0&\hspace{5mm}  e_{1}e_{3}=0  &\hspace{5mm}  e_{2}e_{3}=\alpha e_{2}
\end{array}$$
with $\alpha\neq 0$. Let us denote by $B(\alpha)$ the algebra
defined by means of \textbf{Table $TI_{11}$}. Consequently,
algebra $B(\alpha)$ is isomorphic to an algebra of type
\textbf{A4}, i.e. $B(\alpha)\cong A_{1}(0,\alpha)\cong
A_{1}(\alpha,0)$ for $\alpha\notin\{0,\frac{1}{2}\}$. Similarly,
algebra $B(\frac{1}{2})$ is isomorphic to an algebra of type
\textbf{A2}, i.e. $B(\frac{1}{2})\cong A_{1}(0,\frac{1}{2})$.

\vspace{5mm} \emph{Case} $I_{12}\ \  \delta=-\alpha$

\vspace{3mm} By taking basis $(\gamma e_{1}-\alpha e_{2},\ e_{2},\
e_{3})$ the multiplication table for $A$ becomes

\vspace{3mm}

$$\begin{array}{llll}
  \textbf{Table $TI_{12}$}\hspace{6mm} & \hspace{5mm} e_{1}^{2}=0  &\hspace{5mm} e_{2}^{2}=0&\hspace{5mm} e_{3}^{2}=e_{3} \\
   &\hspace{5mm}  e_{1}e_{2}=0&\hspace{5mm}  e_{1}e_{3}=0 &\hspace{5mm}  e_{2}e_{3}= e_{1}
\end{array}$$

\vspace{5mm}\noindent \textbf{A22)}\ \ \emph{Properties of algebra
$A_{22}$}

\vspace{3mm} $\bullet$ $Ann\ A=\mathbb{R}e_{1},\
\mathcal{N}(A)=Span_{\mathbb{R}}\{e_{1}, e_{2}\},
\mathcal{I}(A)=\{e_{3}\}$,

$\bullet$ 1-dimensional subalgebras: $\mathbb{R}u$ for $u\in
\mathcal{N}(A)\cup \mathcal{I}(A)$,

$\bullet$ 2-dimensional subalgebras: $Span_{\mathbb{R}}\{e_{1},
e_{2}\} , \ Span_{\mathbb{R}}\{e_{1}, e_{3}\}$,

$\bullet$ ideals: $\mathbb{R}e_{1}$, $Span_{\mathbb{R}}\{e_{1},
e_{2}\},\ Span_{\mathbb{R}}\{e_{1}, e_{3}\}$,

$\bullet$ $A^{2}=Span_{\mathbb{R}}\{e_{1}, e_{3}\}$; $A/A^{2}$ is
a 1-dimensional non-null algebra,

$\bullet$ $Der\ A=\left\{
\left[\begin{array}{lll}x&y&0\\0&x&0\\0&0&0\end{array}\right]\ |\
x,y\in \mathbb{R}\right\}$,

$\bullet$ $Aut\ A=\left\{
\left[\begin{array}{lll}x&y&0\\0&x&0\\0&0&1\end{array}\right]\ |\
x,y\in \mathbb{R},\ x\neq 0\right\}$,

$\bullet$ the partition $\mathcal{P}_{A}$ of $\mathbb{R}^{3}$,
defined by the lattice of subalgebras of $A$, consists of:

\vspace{1mm} \hspace{5mm} $\diamond$  the singletons covering the
plane $x^{1}Ox^{2}$,

\vspace{1mm} \hspace{5mm} $\diamond$  the half-axes of axis
$Ox^{3}$ delimited by $O$,

\vspace{1mm} \hspace{5mm} $\diamond$ the quarters of plane
delimited by axes $Ox^{1}$ and $Ox^{3}$ on the plane
$x^{1}Ox^{3}$,

\vspace{1mm} \hspace{5mm} $\diamond$ the quarters of space
$\mathbb{R}^{3}$ delimited by planes $x^{1}Ox^{2}$ and
$x^{1}Ox^{3}$,

\vspace{1mm}$\bullet$ the partition $\mathcal{P}_{A}$ of $A$
induces a partition on the set of integral curves of the
associated homogeneous quadratic differential system (HQDS)
consisting of:

\vspace{1mm} \hspace{5mm} $\diamond$  the singletons consisting of
singular solutions that cover the plane $x^{1}Ox^{2}$,

\vspace{1mm} \hspace{5mm} $\diamond$  the families of
ray-solutions lying on the half-axes of $Ox^{3}$ delimited by $O$,

\vspace{1mm} \hspace{5mm} $\diamond$ the integral curves contained
in the quarters of plane $x^{1}Ox^{3}$ delimited by axes $Ox^{1}$
and $Ox^{3}$,

\vspace{1mm} \hspace{5mm} $\diamond$  the integral curves
contained in the quarters of of space $\mathbb{R}^{3}$ delimited
by planes $x^{1}Ox^{2}$ and $x^{1}Ox^{3}$.

\vspace{3mm} Ideal $Span_{\mathbb{R}}\{e_{1}, e_{3}\}$ compels
each nonsingular integral curve to be torsion-free.

\vspace{3mm} \emph{Subcase $I_{2}$} $\alpha\delta-\beta\gamma= 0,\
\ \alpha\beta\gamma\delta= 0$

\vspace{3mm} We have to consider the cases:

\vspace{2mm}1) $\alpha\neq 0,\ \beta=0,\ \gamma=0,\ \delta=0$,

2) $\alpha\neq 0,\ \beta=0,\ \gamma\neq 0,\ \delta=0$,

3) $\alpha =0,\ \beta\neq 0,\ \gamma=0,\ \delta=0$,

4) $\alpha\neq 0,\ \beta\neq 0,\ \gamma=0,\ \delta=0$,

5) $\alpha =0,\ \beta=0,\ \gamma=0,\ \delta=0$,

6) $\alpha =0,\ \beta=0,\ \gamma=0,\ \delta\neq 0$,

7) $\alpha= 0,\ \beta=0,\ \gamma\neq 0,\ \delta=0$,

8) $\alpha= 0,\ \beta=0,\ \gamma\neq 0,\ \delta\neq 0$.

\begin{proposition}The algebras $1), 2), 4), 6), 7), 8)$ are isomorphic to $\textbf{A11}$,
while algebras $3)$ and $5)$ are respectively isomorphic to
algebras $\textbf{A22}$ and $\textbf{A8}$.
\end{proposition}

\vspace{3mm} \emph{Subcase $II$:}\ \ $\alpha\delta-\beta\gamma\neq
0$

\vspace{3mm} We shall look for a new basis such that the
corresponding multiplication table of algebra should contain a
minimal number of parameters. This is in fact the same with
solving the problem: is or is not semisimple the endomorphism
$L_{e_{3}}$? In fact, this problem is equivalent to the problem:
is or is not diagonalisable the matrix
$$P=\left[\begin{array}{ll}\alpha&\gamma\\ \beta&\delta \end{array} \right]?$$

We have to analyze the next two mutually exclusive cases:

\vspace{2mm}\ (i) $P$ has complex eigenvalues (i.e.
$(\alpha+\delta)^{2}-4(\alpha\delta-\beta\gamma)<0$),

(ii) $P$ has real eigenvalues (i.e.
$(\alpha+\delta)^{2}-4(\alpha\delta-\beta\gamma)\geq 0$).

\vspace{5mm}\emph{Case} (i)

\vspace{3mm}\noindent There exists a basis such that the
multiplication table of algebra has the form:
$$\begin{array}{llll}
  \textbf{Table $\textbf{T5}$}\hspace{6mm} & \hspace{5mm} e_{1}^{2}=0  &\hspace{5mm} e_{2}^{2}=0&\hspace{5mm} e_{3}^{2}=e_{3} \\
   &\hspace{5mm}  e_{1}e_{2}=0&\hspace{5mm}  e_{1}e_{3}=ae_{1}-be_{2} &\hspace{5mm}  e_{2}e_{3}= be_{1}+ae_{2}
\end{array}$$
with $b>0$. Let us denote by $A_{5}(a,b)$ any algebra having
multiplication table $\textbf{T5}$.

\begin{proposition}The algebras $A_{5}(a,b)$ ($b>0$) and $A_{5}(a',b')$ ($b'>0$) are isomorphic if and only if $a=a',\ b=b'$. \end{proposition}
\emph{Proof.} Suppose that the multiplication table of
$A_{5}(a',b')$ is
$$\begin{array}{llll}
  \textbf{Table $\textbf{T'5}$}\hspace{6mm} & \hspace{5mm} f_{1}^{2}=0  &\hspace{5mm} f_{2}^{2}=0&\hspace{5mm} f_{3}^{2}=f_{3} \\
   &\hspace{5mm}  f_{1}f_{2}=0&\hspace{5mm}  f_{1}f_{3}=a'f_{1}-b'f_{2} &\hspace{5mm}  f_{2}f_{3}= b'f_{1}+a'f_{2}
\end{array}$$
and $T(e_{j})=\sum_{i=1}^{3}s_{ij}f_{i}$ is an automorphism of
$A(a,b)$ ($b>0$) with $A(a',b')$ ($b'>0$). Then, $s_{31}=
s_{32}=0,\ s_{33}=1$ and $$s_{13}=2a's_{13}+2b's_{23},\
s_{23}=2a's_{23}-2b's_{13}\Leftrightarrow s_{13}=s_{23}=0.$$ The
conditions $T(e_{1}e_{3})=T(e_{1})T(e_{3})$ and
$T(e_{2}e_{3})=T(e_{2})T(e_{3})$ are equivalent to next equations
in unknown entries $s_{11},\ s_{21},\ s_{12},\ s_{22}$:
$$\left\{ \begin{array}{l}(a-a')s_{11}-b's_{21}-bs_{12}=0\\
b's_{11}+(a-a')s_{21}-bs_{22}=0\\
bs_{11}+(a-a')s_{12}-b's_{22}=0\\
bs_{21}+b's_{12}+(a-a')s_{22}=0. \end{array} \right.$$ This system
has nonzero solution if and only if $a=a'$ and $b=b'$. Indeed, the
determinant
$(b^{2}-b'^{2})^{2}+(a-a')^{2}[(a-a')^{2}+2b^{2}+2b'^{2}]$ of
matrix of coefficients of this homogeneous system is zero if and
only if $a=a'$ and $b=b'$. $\hfill\Box$

\vspace{5mm}\noindent \textbf{A23)}  \emph{Properties of algebra
$A$ of type $A_{5}(a,b)$} \ ($b>0$)

\vspace{3mm} $\bullet$ $Ann\ A=\{0\},\
\mathcal{N}(A)=Span_{\mathbb{R}}\{e_{1}, e_{2}\},
\mathcal{I}(A)=\{e_{3}\}$,

$\bullet$ 1-dimensional subalgebras: $\mathbb{R}u$ for $u\in
\mathcal{N}(A)\cup \mathcal{I}(A)$,

$\bullet$ 2-dimensional subalgebras: $Span_{\mathbb{R}}\{e_{1},
e_{2}\}$ ,

$\bullet$ ideals:  $Span_{\mathbb{R}}\{e_{1}, e_{2}\}$,

$\bullet$ $A^{2}=A$

$\bullet$ $Der\ A=\left\{
\left[\begin{array}{rll}x&y&0\\-y&x&0\\0&0&0\end{array}\right]\ |\
x,y\in \mathbb{R}\right\}\cong \mathbb{C}(\cdot)$,

$\bullet$ $Aut\ A=\left\{
\left[\begin{array}{rll}x&y&0\\-y&x&0\\0&0&1\end{array}\right]\ |\
x,y\in \mathbb{R}^{\ast}\right\}\cong \mathbb{C}^{\ast}(\cdot)$,

$\bullet$ the partition $\mathcal{P}_{A}$ of $\mathbb{R}^{3}$,
defined by the lattice of subalgebras of $A$, consists of:

\vspace{1mm} \hspace{5mm} $\diamond$  the singletons covering the
plane $x^{1}Ox^{2}$,

\vspace{1mm} \hspace{5mm} $\diamond$  the half-axes of axis
$Ox^{3}$ delimited by $O$,

\vspace{1mm} \hspace{5mm} $\diamond$ the half-spaces delimited by
plane $x^{1}Ox^{2}$ less the points of axis $Ox^{3}$,

\vspace{1mm}$\bullet$ the partition $\mathcal{P}_{A}$ of $A$
induces a partition on the set of integral curves of the
associated homogeneous quadratic differential system (HQDS)
consisting of:

\vspace{1mm} \hspace{5mm} $\diamond$  the singletons consisting of
singular solutions that cover the plane $x^{1}Ox^{2}$,

\vspace{1mm} \hspace{5mm} $\diamond$  the families of
ray-solutions lying on the half-axes of $Ox^{3}$ delimited by $O$,

\vspace{1mm} \hspace{5mm} $\diamond$ the integral curves contained
in half-spaces delimited by plane $x^{1}Ox^{2}$ less the points of
axis $Ox^{3}$.

\vspace{5mm}Case (ii)

\vspace{3mm}There exists a basis, consisting of the eigenvectors
of endomorphism $L_{e_{3}}$, such that the multiplication table of
algebra has the form:

$$\begin{array}{llll}
  \textbf{Table\ Tii} \hspace{6mm} & \hspace{5mm} e_{1}^{2}=0  &\hspace{5mm} e_{2}^{2}=0&\hspace{5mm} e_{3}^{2}=e_{3} \\
   &\hspace{5mm}  e_{1}e_{2}=0&\hspace{5mm}  e_{1}e_{3}=ae_{1} &\hspace{5mm}  e_{2}e_{3}= be_{2}
\end{array}$$
with $a,b\in \mathbb{R}$. Let us remark that algebras of type
\textbf{Tii} are necessarily algebras of type \textbf{A2}.

\vspace{5mm} \emph{Case Spec\ D=(1,2,0)}

\vspace{3mm} There exists a basis
$\mathcal{B}=(e_{1},e_{2},e_{3})$ such that the multiplication
table of algebra is

$$\begin{array}{llll}
  \textbf{Table T'}\hspace{6mm} & \hspace{5mm} e_{1}^{2}=be_{2}  &\hspace{5mm} e_{2}^{2}=0&\hspace{5mm} e_{3}^{2}=je_{3} \\
   &\hspace{5mm}  e_{1}e_{2}=0  &\hspace{5mm}  e_{1}e_{3}=pe_{1} &\hspace{5mm}  e_{2}e_{3}=te_{2}
\end{array}$$
with $b,j,p,t\in \mathbb{R}$. Will be suitable to use the change
of bases $(e_{1},e_{2},e_{3})\rightarrow (e_{2},e_{1},e_{3})$. The
corresponding multiplication table is:
$$\begin{array}{llll}
  \textbf{Table T"}\hspace{6mm} & \hspace{5mm} e_{1}^{2}=0  &\hspace{5mm} e_{2}^{2}=be_{1}&\hspace{5mm} e_{3}^{2}=je_{3} \\
   &\hspace{5mm}  e_{1}e_{2}=0  &\hspace{5mm}  e_{1}e_{3}=te_{1} &\hspace{5mm}  e_{2}e_{3}=pe_{2}
\end{array}$$
with $b,j,p,t\in \mathbb{R}$.

We have to distinguish two complementary situations:

$$I)\ \ bj\neq 0,\ \ \ II)\ \ bj=0$$

\vspace{3mm} \noindent I)\ \emph{Case} $bj\neq 0$

\vspace{3mm} There exists a basis
$\mathcal{B}=(e_{1},e_{2},e_{3})$ such that the multiplication
table of algebra is
$$\begin{array}{llll}
  \textbf{Table T6}\hspace{6mm} & \hspace{5mm} e_{1}^{2}=0  &\hspace{5mm} e_{2}^{2}=e_{1}&\hspace{5mm} e_{3}^{2}=e_{3} \\
   &\hspace{5mm}  e_{1}e_{2}=0  &\hspace{5mm}  e_{1}e_{3}=\alpha e_{1} &\hspace{5mm}  e_{2}e_{3}=\beta e_{2}
\end{array}$$
with $\alpha,\beta\in \mathbb{R}$. Let us denote by
$A_{6}(\alpha,\beta)$ any algebra having the multiplication table
\textbf{T6}.

\begin{proposition}\label{p16} The algebras $A_{6}(\alpha,\beta)$ and
$A_{6}(\alpha',\beta')$ are isomorphic if and only if
$\alpha=\alpha',\ \beta=\beta'$.
\end{proposition}

\begin{proposition}
Every algebra $A_{6}(\alpha,\beta)$ has:
$$Ann\ A=\left\{\begin{array}{lll}\mathbb{R}e_{1} & if & \alpha=0\\ \{0\}&if&\alpha\neq 0. \end{array}\right.$$
$$\mathcal{N}(A)=\mathbb{R}e_{1}$$

$$\mathcal{I}(A)=\left\{\begin{array}{lll}\{e_{3}\}&if& \alpha\neq \frac{1}{2},\ \beta\neq \frac{1}{2}\\
xe_{1}+e_{3}&if&\alpha=\frac{1}{2}\\
\{\frac{y^{2}}{1-2\alpha}e_{1}+ye_{1}+e_{3}\ |\ y\in
\mathbb{R}\}&if& \alpha\neq \frac{1}{2},\ \beta=\frac{1}{2}.
  \end{array}\right.$$

\end{proposition}

Accordingly, we have to consider the classes of algebras:

\begin{center}(i)\ \ $A_{6}(0,\beta)$ \ and\ \ (ii)\
$A_{6}(\alpha,\beta)\ with\ \alpha\neq 0$\end{center}

\vspace{3mm} \emph{Case (i)}

\vspace{3mm} \emph{Algebras $A_{6}(0,\beta)$ with $\beta\notin
\{\frac{1}{2},0\}$}

\vspace{3mm}Proposition \ref{p16} implies:
\begin{corollary} The algebras $A_{6}(0,\beta)$ and
$A_{6}(0,\beta')$ are isomorphic if and only if $\beta=\beta'$.
\end{corollary}

\vspace{3mm}\noindent \textbf{A24)}  \emph{Properties of algebra
$A=A_{6}(0,\beta)$ with $\beta\notin\{0, \frac{1}{2})$}

\vspace{2mm} $\bullet$ $Ann\ A=\mathbb{R}e_{1},\
\mathcal{N}(A)=\mathbb{R}e_{1},\ \mathcal{I}(A)=\{e_{3}\}$,

$\bullet$ 1-dimensional subalgebras: $\mathbb{R}e_{1},\
\mathbb{R}e_{3}$,

$\bullet$ 2-dimensional subalgebras: $
Span_{\mathbb{R}}\{e_{1},e_{2}\},\
Span_{\mathbb{R}}\{e_{1},e_{3}\}$,

$\bullet$ ideals:  $\mathbb{R}e_{1},\ Span_{\mathbb{R}}\{e_{1},
e_{2}\}$,

$\bullet$ $A^{2}=A$

$\bullet$ $Der\ A=\left\{
\left[\begin{array}{rll}2x&0&0\\0&x&0\\0&0&0\end{array}\right]\ |\
x\in \mathbb{R}\right\}\cong \mathbb{R}D$,

$\bullet$ $Aut\ A=\left\{
\left[\begin{array}{rll}x^{2}&0&0\\0&x&0\\0&0&1\end{array}\right]\
|\ x\in \mathbb{R}^{\ast}\}\right.\cong \mathbb{R}^{\ast}(\cdot)$,

$\bullet$ the partition $\mathcal{P}_{A}$ of $\mathbb{R}^{3}$,
defined by the lattice of subalgebras of $A$, consists of:

\vspace{1mm} \hspace{5mm} $\diamond$  the singletons covering the
axis $Ox^{1}$,

\vspace{1mm} \hspace{5mm} $\diamond$  the half-axes of axis
$Ox^{3}$ delimited by $O$,

\vspace{1mm} \hspace{5mm} $\diamond$ the quarters delimited by
axes $Ox^{1}$ and $Ox^{3}$ on plane $x^{1}Ox^{3}$,

\vspace{1mm} \hspace{5mm} $\diamond$ the half-planes delimited by
axis $Ox^{1}$ on plane $x^{1}Ox^{2}$,

\vspace{1mm} \hspace{5mm} $\diamond$ the quarters of space
delimited by planes  $x^{1}Ox^{2}$ and $x^{1}Ox^{3}$,

\vspace{1mm}$\bullet$ the partition $\mathcal{P}_{A}$ of $A$
induces a partition on the set of integral curves of the
associated homogeneous quadratic differential system (HQDS)
consisting of:

\vspace{1mm} \hspace{5mm} $\diamond$  the singletons consisting of
singular solutions that cover the axis $Ox^{1}$,

\vspace{1mm} \hspace{5mm} $\diamond$  the families of
ray-solutions lying on the half-axes of $Ox^{3}$ delimited by $O$,

\vspace{1mm} \hspace{5mm} $\diamond$ the integral curves contained
in the quarters delimited by axes $Ox^{1}$ and $Ox^{3}$ on plane
$x^{1}Ox^{3}$,

\vspace{1mm} \hspace{5mm} $\diamond$ the integral curves contained
in the half-planes delimited by axis $Ox^{1}$ on plane
$x^{1}Ox^{2}$,

\vspace{1mm} \hspace{5mm} $\diamond$ the integral curves contained
in the quarters of space delimited by planes  $x^{1}Ox^{2}$ and
$x^{1}Ox^{3}$.

\vspace{5mm}\noindent \textbf{A25)}  \emph{Properties of algebra
$A=A_{6}(0,0)$}

\vspace{2mm} $\bullet$ $Ann\ A=\mathbb{R}e_{1},\
\mathcal{N}(A)=\mathbb{R}e_{1},\ \mathcal{I}(A)=\{e_{3}\}$,

$\bullet$ 1-dimensional subalgebras: $\mathbb{R}e_{1},\
\mathbb{R}e_{3}$,

$\bullet$ 2-dimensional subalgebras: $
Span_{\mathbb{R}}\{e_{1},e_{2}\},\
Span_{\mathbb{R}}\{e_{1},e_{3}\}$,

$\bullet$ ideals:  $\mathbb{R}e_{1},\ Span_{\mathbb{R}}\{e_{1},
e_{2}\},\ Span_{\mathbb{R}}\{e_{1},e_{3}\}$,

$\bullet$ $A^{2}=Span_{\mathbb{R}}\{e_{1},e_{3}\}$

$\bullet$ $Der\ A=\left\{
\left[\begin{array}{rll}2x&y&0\\0&x&0\\0&0&0\end{array}\right]\ |\
x,y\in \mathbb{R}\right\}$,

$\bullet$ $Aut\ A=\left\{
\left[\begin{array}{rll}x^{2}&y&0\\0&x&0\\0&0&1\end{array}\right]\
|\ x,y\in \mathbb{R},\ x\neq 0\}\right.$,

$\bullet$ the partition $\mathcal{P}_{A}$ of $\mathbb{R}^{3}$,
defined by the lattice of subalgebras of $A$, consists of:

\vspace{1mm} \hspace{5mm} $\diamond$  the singletons covering the
axis $Ox^{1}$,

\vspace{1mm} \hspace{5mm} $\diamond$  the half-axes of axis
$Ox^{3}$ delimited by $O$,

\vspace{1mm} \hspace{5mm} $\diamond$ the quarters delimited by
axes $Ox^{1}$ and $Ox^{3}$ on plane $x^{1}Ox^{3}$,

\vspace{1mm} \hspace{5mm} $\diamond$ the half-planes delimited by
axis $Ox^{1}$ on plane $x^{1}Ox^{2}$,

\vspace{1mm} \hspace{5mm} $\diamond$ the quarters of space
delimited by planes  $x^{1}Ox^{2}$ and $x^{1}Ox^{3}$,

\vspace{1mm}$\bullet$ the partition $\mathcal{P}_{A}$ of $A$
induces a partition on the set of integral curves of the
associated homogeneous quadratic differential system (HQDS)
consisting of:

\vspace{1mm} \hspace{5mm} $\diamond$  the singletons consisting of
singular solutions that cover the axis $Ox^{1}$,

\vspace{1mm} \hspace{5mm} $\diamond$  the families of
ray-solutions lying on the half-axes of $Ox^{3}$ delimited by $O$,

\vspace{1mm} \hspace{5mm} $\diamond$ the integral curves contained
in the quarters delimited by axes $Ox^{1}$ and $Ox^{3}$ on plane
$x^{1}Ox^{3}$,

\vspace{1mm} \hspace{5mm} $\diamond$ the integral curves contained
in the half-planes delimited by axis $Ox^{1}$ on plane
$x^{1}Ox^{2}$,

\vspace{1mm} \hspace{5mm} $\diamond$ the integral curves contained
in the quarters of space delimited by planes  $x^{1}Ox^{2}$ and
$x^{1}Ox^{3}$.

\vspace{5mm}\noindent \textbf{A26)}  \emph{Properties of algebra
$A=A_{6}(0,\frac{1}{2})$}

\vspace{2mm} $\bullet$ $Ann\ A=\mathbb{R}e_{1},\
\mathcal{N}(A)=\mathbb{R}e_{1},\
\mathcal{I}(A)=\{x^{2}e_{1}+xe_{2}+e_{3}\ |\ x\in \mathbb{R}\}$,

$\bullet$ 1-dimensional subalgebras: $\mathbb{R}u$ for $u\in
\mathcal{N}(A)\cup \mathcal{I}(A)$,

$\bullet$ 2-dimensional subalgebras: $
Span_{\mathbb{R}}\{e_{1},be_{2}+ce_{3}\},\ (b^{2}+c^{2}\neq 0)$,

$\bullet$ ideals:  $\mathbb{R}e_{1},\ Span_{\mathbb{R}}\{e_{1},
e_{2}\}$,

$\bullet$ $A^{2}=A$,

$\bullet$ $Der\ A=\left\{
\left[\begin{array}{rll}2x&2y&0\\0&x&y\\0&0&0\end{array}\right]\
|\ x,y\in \mathbb{R}\right\}$,

$\bullet$ $Aut\ A=\left\{
\left[\begin{array}{rll}x^{2}&2xy&y^{2}\\0&x&y\\0&0&1\end{array}\right]\
|\ x,y\in \mathbb{R},\ x\neq 0\}\right.$,

$\bullet$ the partition $\mathcal{P}_{A}$ of $\mathbb{R}^{3}$,
defined by the lattice of subalgebras of $A$, consists of:

\vspace{1mm} \hspace{5mm} $\diamond$  the singletons covering the
axis $Ox^{1}$,

\vspace{1mm} \hspace{5mm} $\diamond$  the half-lines of
$\mathbb{R}u$ for $u\in \mathcal{I}(A)$ delimited by $O$ (these
lines cover the cone $(x^{2})^{2}-x^{1}x^{3}=0$),

\vspace{1mm} \hspace{5mm} $\diamond$ the connected components
delimited by axis $Ox^{1}$ and the cone $(x^{2})^{2}-x^{1}x^{3}=0$
on each plane containing $Ox^{1}$ (of course, this cone has a real
contribution whenever $x^{1}x^{3}>0$),

\vspace{1mm}$\bullet$ the partition $\mathcal{P}_{A}$ of $A$
induces a partition on the set of integral curves of the
associated homogeneous quadratic differential system (HQDS)
consisting of:

\vspace{1mm} \hspace{5mm} $\diamond$  the singletons consisting of
singular solutions that cover the axis $Ox^{1}$,

\vspace{1mm} \hspace{5mm} $\diamond$  the families of
ray-solutions lying on the half-lines of $\mathbb{R}u$ for $u\in
\mathcal{I}(A)$ delimited by $O$,

\vspace{1mm} \hspace{5mm} $\diamond$ the integral curves contained
in the connected components delimited by axis $Ox^{1}$ and the
cone $(x^{2})^{2}-x^{1}x^{3}=0$ on each plane containing $Ox^{1}$.

\vspace{3mm} \emph{Case (ii)}

\vspace{3mm} \emph{Algebras $A_{6}(\alpha,\beta)$ with $\alpha\neq
0$}

\vspace{3mm}Proposition \ref{p16} implies:
\begin{corollary} The algebras $A_{6}(\alpha,\beta)\ (\alpha\neq
0)$\ and\ $A_{6}(\alpha ',\beta ')\ (\alpha'\neq 0)$ are
isomorphic if and only if\ $\alpha=\alpha ',\ \beta=\beta '$.
\end{corollary}

\begin{proposition}
Every algebra $A_{6}(\alpha,\beta)$ has:
$$Ann\ A=\{0\}$$
$$\mathcal{N}(A)=\mathbb{R}e_{1}$$
$$\mathcal{I}(A)=\left\{\begin{array}{lll}\{e_{3}\}&if& \alpha\neq \frac{1}{2},\ \beta\neq \frac{1}{2}\\
\{xe_{1}+e_{3}\ |\ x\in \mathbb{R}\}&if& \alpha=\frac{1}{2},\
\beta\neq  \frac{1}{2}\\
\{\frac{y^2}{1-2\alpha}e_{1}+ye_{2}+e_{3}\ |\ y\in
\mathbb{R}\}&if& \alpha \neq \frac{1}{2},\
\beta =\frac{1}{2}\\
\{xe_{1}+e_{3}\ |\ x\in \mathbb{R}\}&if& \alpha = \frac{1}{2},\
\beta= \frac{1}{2}.\end{array} \right.$$
\end{proposition}

\vspace{3mm}\noindent \textbf{A27)} \emph{Properties of algebras
$A$ of type $A_{6}(\alpha,\beta)$ with $\alpha,\ \beta\notin
\{0,\frac{1}{2}\}$ and $\alpha\neq \beta$}

\vspace{2mm} $\bullet$ $Ann\ A=\{0\},\
\mathcal{N}(A)=\mathbb{R}e_{1},\ \mathcal{I}(A)=\{e_{3}\}$,

$\bullet$ 1-dimensional subalgebras: $\mathbb{R}e_{1},\
\mathbb{R}e_{3}$ ,

$\bullet$ 2-dimensional subalgebras:
$Span_{\mathbb{R}}\{e_{1},e_{2}\},\
Span_{\mathbb{R}}\{e_{1},e_{3}\}$,

$\bullet$ ideals:  $\mathbb{R}e_{1},\ Span_{\mathbb{R}}\{e_{1},
e_{2}\}$,

$\bullet$ $A^{2}=A$

$\bullet$ $Der\ A=\left\{
\left[\begin{array}{rll}2x&0&0\\0&x&0\\0&0&0\end{array}\right]\ |\
x\in \mathbb{R}\right\}\cong \mathbb{R}D$,

$\bullet$ $Aut\ A=\left\{
\left[\begin{array}{rll}x^{2}&0&0\\0&x&0\\0&0&1\end{array}\right]\
|\ x\in \mathbb{R}^{\ast}\right\}\cong\mathbb{R}^{\ast}(\cdot)$,

$\bullet$ the partition $\mathcal{P}_{A}$ of $\mathbb{R}^{3}$,
defined by the lattice of subalgebras of $A$, consists of:

\vspace{1mm} \hspace{5mm} $\diamond$  the singletons covering the
axis $Ox^{1}$,

\vspace{1mm} \hspace{5mm} $\diamond$  the half-axes of $Ox^{3}$
delimited by $O$,

\vspace{1mm} \hspace{5mm} $\diamond$ the half-planes delimited by
axis $Ox^{1}$ on plane $x^{1}Ox^{2}$,

\vspace{1mm} \hspace{5mm} $\diamond$ the quarters of plane
delimited by axes $Ox^{1}$ and $Ox^{3}$ on plane $x^{1}Ox^{3}$,

\vspace{1mm} \hspace{5mm} $\diamond$ the quarters of space
delimited by planes $x^{1}Ox^{2}$ and $x^{1}Ox^{3}$,

\vspace{1mm}$\bullet$ the partition $\mathcal{P}_{A}$ of $A$
induces a partition on the set of integral curves of the
associated homogeneous quadratic differential system (HQDS)
consisting of:

\vspace{1mm} \hspace{5mm} $\diamond$  the singletons consisting of
singular solutions that cover the axis $Ox^{1}$,

\vspace{1mm} \hspace{5mm} $\diamond$  the families of
ray-solutions lying on the half-axes of $Ox^{3}$ delimited by $O$,

\vspace{1mm} \hspace{5mm} $\diamond$ the integral curves contained
in the half-planes delimited by axis $Ox^{1}$ on plane
$x^{1}Ox^{2}$,

\vspace{1mm} \hspace{5mm} $\diamond$ the integral curves contained
in the quarters of plane delimited by axes $Ox^{1}$ and $Ox^{3}$
on plane $x^{1}Ox^{3}$,

\vspace{1mm} \hspace{5mm} $\diamond$ the integral curves contained
in the quarters of space delimited by planes $x^{1}Ox^{2}$ and
$x^{1}Ox^{3}$.

\vspace{5mm}\noindent \textbf{A28)} \emph{Properties of algebras
$A$ of type $A_{6}(\alpha,\alpha)$ with $\alpha\notin\{0,\
\frac{1}{2}\}$}

\vspace{3mm}$\bullet$ the algebras $A_{6}(\alpha,\alpha)\
(\alpha\neq 0)$\ and\ $A_{6}(\alpha ',\alpha ')\ (\alpha'\neq 0)$
are isomorphic if and only if\ $\alpha=\alpha '$

$\bullet$ $Ann\ A=\{0\},\ \mathcal{N}(A)=\mathbb{R}e_{1},\
\mathcal{I}(A)=\{e_{3}\}$,

$\bullet$ 1-dimensional subalgebras: $\mathbb{R}e_{1},\
\mathbb{R}e_{3}$,

$\bullet$ 2-dimensional subalgebras:
$Span_{\mathbb{R}}\{e_{1},e_{2}\},\
Span_{\mathbb{R}}\{e_{1},e_{3}\}$,

$\bullet$ ideals:  $\mathbb{R}e_{1},\ Span_{\mathbb{R}}\{e_{1},
e_{2}\}$,

$\bullet$ $A^{2}=A$

$\bullet$ $Der\ A=\left\{
\left[\begin{array}{rll}2x&0&0\\y&x&0\\0&0&0\end{array}\right]\ |\
x,y\in \mathbb{R}\right\}$,

$\bullet$ $Aut\ A=\left\{
\left[\begin{array}{rll}x^{2}&0&0\\y&x&0\\0&0&1\end{array}\right]\
|\ x,y\in \mathbb{R},\ x\neq 0\}\right.$,

$\bullet$ the partition $\mathcal{P}_{A}$ of $\mathbb{R}^{3}$,
defined by the lattice of subalgebras of $A$, consists of:

\vspace{1mm} \hspace{5mm} $\diamond$  the singletons covering the
axis $Ox^{1}$,

\vspace{1mm} \hspace{5mm} $\diamond$  the half-axes of $Ox^{3}$
delimited by $O$,

\vspace{1mm} \hspace{5mm} $\diamond$ the half-planes delimited by
axis $Ox^{1}$ on plane $x^{1}Ox^{2}$,

\vspace{1mm} \hspace{5mm} $\diamond$ the quarters of plane
delimited by axes $Ox^{1}$ and $Ox^{3}$ on plane $x^{1}Ox^{3}$,

\vspace{1mm} \hspace{5mm} $\diamond$ the quarters of space
delimited by planes $x^{1}Ox^{2}$ and $x^{1}Ox^{3}$,

\vspace{1mm}$\bullet$ the partition $\mathcal{P}_{A}$ of $A$
induces a partition on the set of integral curves of the
associated homogeneous quadratic differential system (HQDS)
consisting of:

\vspace{1mm} \hspace{5mm} $\diamond$  the singletons consisting of
singular solutions that cover the axis $Ox^{1}$,

\vspace{1mm} \hspace{5mm} $\diamond$  the families of
ray-solutions lying on the half-axes of $Ox^{3}$ delimited by $O$,

\vspace{1mm} \hspace{5mm} $\diamond$ the integral curves contained
in the half-planes delimited by axis $Ox^{1}$ on plane
$x^{1}Ox^{2}$,

\vspace{1mm} \hspace{5mm} $\diamond$ the integral curves contained
in the quarters of plane delimited by axes $Ox^{1}$ and $Ox^{3}$
on plane $x^{1}Ox^{3}$,

\vspace{1mm} \hspace{5mm} $\diamond$ the integral curves contained
in the quarters of space delimited by planes $x^{1}Ox^{2}$ and
$x^{1}Ox^{3}$.

\vspace{5mm}\noindent \textbf{A29)} \emph{Properties of algebras
$A$ of type $A_{6}(\alpha,\frac{1}{2})$ with} $\alpha\notin\{0,
\frac{1}{2}\}$

\vspace{3mm}$\bullet$ the algebras $A_{6}(\alpha,\frac{1}{2})\
(\alpha\notin\{0, \frac{1}{2}\})$\ and\ $A_{6}(\alpha
',\frac{1}{2})\ (\alpha'\notin\{0, \frac{1}{2}\})$ are isomorphic
if and only if\ $\alpha=\alpha '$,

$\bullet$ $Ann\ A=\{0\},\ \mathcal{N}(A)=\mathbb{R}e_{1},\
\mathcal{I}(A)=\{\frac{y^{2}}{1-2\alpha}e_{1}+ye_{2}+e_{3}\ |\
y\in \mathbb{R}\}$,

$\bullet$ 1-dimensional subalgebras: $\mathbb{R}u$ for $u\in
\mathcal{N}(A)\cup \mathcal{I}(A)$,

$\bullet$ 2-dimensional subalgebras:
$Span_{\mathbb{R}}\{e_{1},be_{2}+ce_{3}\},\ b^2+c^2\neq 0,$

$\bullet$ ideals:  $\mathbb{R}e_{1},\ Span_{\mathbb{R}}\{e_{1},
e_{2}\}$,

$\bullet$ $A^{2}=A$

$\bullet$ $Der\ A=\left\{
\left[\begin{array}{rll}2x&0&0\\0&x&y\\0&0&0\end{array}\right]\ |\
x,y\in \mathbb{R}\right\}$,

$\bullet$ $Aut\ A=\left\{
\left[\begin{array}{rll}x^{2}&0&0\\0&x&y\\0&0&1\end{array}\right]\
|\ x,y\in \mathbb{R},\ x\neq 0\}\right.$,

$\bullet$ the partition $\mathcal{P}_{A}$ of $\mathbb{R}^{3}$,
defined by the lattice of subalgebras of $A$, consists of:

\vspace{1mm} \hspace{5mm} $\diamond$  the singletons covering the
axis $Ox^{1}$,

\vspace{1mm} \hspace{5mm} $\diamond$ the half-lines of lines
$\mathbb{R}u$ for $u\in \mathcal{I}(A)$ delimited by $O$ (these
lines cover the cone $(x^{2})^{2}-(1-2\alpha)x^{1}x^{3}=0$),

\vspace{1mm} \hspace{5mm} $\diamond$ the connected components
delimited by axis $Ox^{1}$ and the cone
$(x^{2})^{2}-(1-2\alpha)x^{1}x^{3}=0$ on each plane containing
$Ox^{1}$,

\vspace{1mm}$\bullet$ the partition $\mathcal{P}_{A}$ of $A$
induces a partition on the set of integral curves of the
associated homogeneous quadratic differential system (HQDS)
consisting of:

\vspace{1mm} \hspace{5mm} $\diamond$  the singletons consisting of
singular solutions that cover the axis $Ox^{1}$,

\vspace{1mm} \hspace{5mm} $\diamond$  the families of
ray-solutions lying on the half-lines of lines $\mathbb{R}u$ for
$u\in \mathcal{I}(A)$ delimited by $O$ (these lines cover the cone
$(x^{2})^{2}-(1-2\alpha)x^{1}x^{3}=0$).

\vspace{1mm} \hspace{5mm} $\diamond$ the integral curves contained
in the connected components delimited by axis $Ox^{1}$ and the
cone $(x^{2})^{2}-(1-2\alpha)x^{1}x^{3}=0$ on each plane
containing $Ox^{1}$.

\vspace{5mm}\noindent \textbf{A30)} \emph{Properties of algebras
$A$ of type $A_{6}(\frac{1}{2},\frac{1}{2})$}

\vspace{2mm} $\bullet$ $Ann\ A=\{0\},\
\mathcal{N}(A)=\mathbb{R}e_{1},\ \mathcal{I}(A)=\{xe_{1}+e_{3}\ |\
x\in \mathbb{R}\}$,

$\bullet$ 1-dimensional subalgebras: $\mathbb{R}u$ for $u\in
\mathcal{N}(A)\cup \mathcal{I}(A)$,

$\bullet$ 2-dimensional subalgebras:
$Span_{\mathbb{R}}\{e_{1},be_{2}+ce_{3}\}\ (b^{2}+c^{2}\neq 0) $,

$\bullet$ ideals:  $\mathbb{R}e_{1},\ Span_{\mathbb{R}}\{e_{1},
e_{2}\}$,

$\bullet$ $A^{2}=A$

$\bullet$ $Der\ A=\left\{
\left[\begin{array}{rll}2x&0&0\\y&x&z\\0&0&0\end{array}\right]\ |\
x,y,z\in \mathbb{R}\right\}$,

$\bullet$ $Aut\ A=\left\{
\left[\begin{array}{rll}x^{2}&0&0\\y&x&z\\0&0&1\end{array}\right]\
|\ x,y,z\in \mathbb{R},\ x\neq 0\}\right.$,

$\bullet$ the partition $\mathcal{P}_{A}$ of $\mathbb{R}^{3}$,
defined by the lattice of subalgebras of $A$, consists of:

\vspace{1mm} \hspace{5mm} $\diamond$  the singletons covering the
axis $Ox^{1}$,

\vspace{1mm} \hspace{5mm} $\diamond$ the half-lines of lines
$\mathbb{R}u$ for $u\in \mathcal{I}(A)$ delimited by $O$ (these
lines cover the plane $x^{1}Ox^{3}$ less axis $Ox^{1}$),

\vspace{1mm} \hspace{5mm} $\diamond$ the half-planes delimited by
axis $Ox^{1}$ on each plane passing through $Ox^{1}$ without
$x^{1}Ox^{3}$,

\vspace{1mm}$\bullet$ the partition $\mathcal{P}_{A}$ of $A$
induces a partition on the set of integral curves of the
associated homogeneous quadratic differential system (HQDS)
consisting of:

\vspace{1mm} \hspace{5mm} $\diamond$  the singletons consisting of
singular solutions that cover the axis $Ox^{1}$,

\vspace{1mm} \hspace{5mm} $\diamond$  the families of
ray-solutions lying on the half-lines of lines $\mathbb{R}u$ for
$u\in \mathcal{I}(A)$ delimited by $O$ (these lines cover the
plane $x^{1}Ox^{3}$ less axis $Ox^{1}$),

\vspace{1mm} \hspace{5mm} $\diamond$ the integral curves contained
in the half-planes delimited by axis $Ox^{1}$ on  each plane
passing through $Ox^{1}$ without $x^{1}Ox^{3}$.

\vspace{2mm} Note that each nonsingular integral curve has zero
torsion.

\vspace{5mm}\noindent \textbf{31)} \emph{Properties of algebras
$A$ of type $A_{6}(\frac{1}{2},\beta)$ with
$\beta\notin\{0,\frac{1}{2}\}$}

\vspace{3mm}$\bullet$ the algebras $A_{6}(\frac{1}{2},\beta)\
(\beta\neq \frac{1}{2})$\ and\ $A_{6}(\frac{1}{2},\beta')\
(\beta'\neq \frac{1}{2})$ are isomorphic if and only if\
$\beta=\beta '$,

$\bullet$ $Ann\ A=\{0\},\ \mathcal{N}(A)=\mathbb{R}e_{1},\
\mathcal{I}(A)=\{xe_{1}+e_{3}\ |\ x\in \mathbb{R}\}$,

$\bullet$ 1-dimensional subalgebras: $\mathbb{R}u$ for $u\in
\mathcal{N}(A)\cup \mathcal{I}(A)$,

$\bullet$ 2-dimensional subalgebras:
$Span_{\mathbb{R}}\{e_{1},e_{2}\},\
Span_{\mathbb{R}}\{e_{1},e_{3}\} $,

$\bullet$ ideals:  $\mathbb{R}e_{1},\ Span_{\mathbb{R}}\{e_{1},
e_{2}\}$,

$\bullet$ $A^{2}=A$

$\bullet$ $Der\ A=\left\{
\left[\begin{array}{rcc}2x&0&(\frac{1}{2}-\beta)y\\y&x&0\\0&0&0\end{array}\right]\
|\ x,y\in \mathbb{R}\right\}$,

$\bullet$ $Aut\ A=\left\{
\left[\begin{array}{rcc}x^{2}&0&y\\\frac{2xy}{1-2\beta}&x&\frac{y^{2}}{1-2\beta}\\0&0&1\end{array}\right]\
|\ x,y\in \mathbb{R},\ x\neq 0\}\right.$,

$\bullet$ the partition $\mathcal{P}_{A}$ of $\mathbb{R}^{3}$,
defined by the lattice of subalgebras of $A$, consists of:

\vspace{1mm} \hspace{5mm} $\diamond$  the singletons covering the
axis $Ox^{1}$,

\vspace{1mm} \hspace{5mm} $\diamond$ the half-lines of lines
$\mathbb{R}u$ for $u\in \mathcal{I}(A)$ delimited by $O$ (these
lines cover the plane $x^{1}Ox^{3}$ less the axis $Ox^{1}$),

\vspace{1mm} \hspace{5mm} $\diamond$ the half-planes delimited by
axis $Ox^{1}$ on the plane $x^{1}Ox^{2}$,

\vspace{1mm} \hspace{5mm} $\diamond$ the connected components of
space delimited by planes $x^{1}Ox^{2}$ and $x^{1}Ox^{3}$,

\vspace{1mm}$\bullet$ the partition $\mathcal{P}_{A}$ of $A$
induces a partition on the set of integral curves of the
associated homogeneous quadratic differential system (HQDS)
consisting of:

\vspace{1mm} \hspace{5mm} $\diamond$  the singletons consisting of
singular solutions that cover the axis $Ox^{1}$,

\vspace{1mm} \hspace{5mm} $\diamond$  the families of
ray-solutions lying on the half-lines of lines $\mathbb{R}u$ for
$u\in \mathcal{I}(A)$ delimited by $O$ (these lines cover the
plane $x^{1}Ox^{3}$ less the axis $Ox^{1}$),

\vspace{1mm} \hspace{5mm} $\diamond$ the integral curves contained
in the connected components of space delimited by planes
$x^{1}Ox^{2}$ and $x^{1}Ox^{3}$.

\vspace{2mm} Note that each nonsingular integral curve has zero
torsion.

\vspace{5mm} \emph{Case II: \ \ bj=0}

\vspace{3mm} We have to consider the subcases:
$$(i)\ \ b=0,\ j\neq 0,\ \ \ (ii)\ \ b\neq 0,\ j=0,\ \ \ (iii)\ \ b=j=0.$$

\vspace{3mm} \noindent\emph{Subcase (i)}

\vspace{3mm} There exists a basis
$\mathcal{B}=(e_{1},e_{2},e_{3})$ such that the multiplication
table of algebra becomes

$$\begin{array}{lll}
   \hspace{5mm} e_{1}^{2}=0  &\hspace{5mm} e_{2}^{2}=0&\hspace{5mm} e_{3}^{2}=e_{3} \\
    \hspace{5mm} e_{1}e_{2}=0&\hspace{5mm}  e_{1}e_{3}=\alpha e_{1} &\hspace{5mm}  e_{2}e_{3}=\beta e_{2}
\end{array}$$
with $\alpha,\beta\in \mathbb{R}$. Consequently, this algebra is
isomorphic to an algebra of type $\textbf{A14}$.

\vspace{3mm} \noindent\emph{Subcase (ii)}

\vspace{3mm} There exists a basis
$\mathcal{B}=(e_{1},e_{2},e_{3})$ such that the multiplication
table of algebra has form

$$\begin{array}{llll}
  \textbf{Table T7}\hspace{6mm} & \hspace{5mm} e_{1}^{2}=e_{2}  &\hspace{5mm} e_{2}^{2}=e_{1}&\hspace{5mm} e_{3}^{2}=0 \\
   &\hspace{5mm}  e_{1}e_{2}=0  &\hspace{5mm}  e_{1}e_{3}=\alpha e_{1} &\hspace{5mm}  e_{2}e_{3}=\beta e_{2}
\end{array}$$
with $\alpha, \beta \in \mathbb{R}$. The algebra having the
multiplication table $\textbf{T7}$ is denoted by $A_{7}(\alpha,
\beta)$.

It is proved that only the following three multiplication tables
are of interest:

$$\begin{array}{llll}
  (ii_{1})& \hspace{5mm} e_{1}^{2}=e_{2}  &\hspace{5mm} e_{2}^{2}=0&\hspace{5mm} e_{3}^{2}=0 \\
  \hspace{5mm} & \hspace{5mm} e_{1}e_{2}=0&\hspace{5mm}  e_{1}e_{3}=0 &\hspace{5mm}  e_{2}e_{3}=0\\
  \\
(ii_{2})& \hspace{5mm} e_{1}^{2}=e_{2}  &\hspace{5mm} e_{2}^{2}=0&\hspace{5mm} e_{3}^{2}=0 \\
  \hspace{5mm} & \hspace{5mm} e_{1}e_{2}=0&\hspace{5mm}  e_{1}e_{3}=e_{1} &\hspace{5mm}  e_{2}e_{3}=\beta e_{2}\\
  \\

(ii_{3})& \hspace{5mm} e_{1}^{2}=e_{2}  &\hspace{5mm} e_{2}^{2}=0&\hspace{5mm} e_{3}^{2}=0 \\
  \hspace{5mm} & \hspace{5mm} e_{1}e_{2}=0&\hspace{5mm}  e_{1}e_{3}=\alpha e_{1} &\hspace{5mm}
  e_{2}e_{3}=e_{2}.
\end{array}$$

\vspace{5mm}\noindent \textbf{A32)} \emph{Properties of algebra
$A=A_{7}(0,0)$ (of type $(ii_{1})$)}

\vspace{2mm} $\bullet$ $Ann\ A=Span_{\mathbb{R}}\{e_{2},e_{3}\},\
\mathcal{N}(A)=Span_{\mathbb{R}}\{e_{2},e_{3}\},\
\mathcal{I}(A)=\emptyset$,

$\bullet$ 1-dimensional subalgebras: $\mathbb{R}u$ for $u\in
\mathcal{N}(A)$,

$\bullet$ 2-dimensional subalgebras:
$Span_{\mathbb{R}}\{e_{2},ae_{1}+ce_{3}\}\ (a^{2}+c^{2}\neq 0) $,

$\bullet$ ideals: $\mathbb{R}u$ for $u\in Ann\ A$,
$Span_{\mathbb{R}}\{e_{2},ae_{1}+ce_{3}\}\ (a^{2}+c^{2}\neq 0) $,

$\bullet$ $A^{2}=\mathbb{R}e_{2}$; $A$ is nilpotent; $A/A^{2}$ is
a 2-dimensional null algebra,

$\bullet$ $A$ is a nilpotent commutative associative algebra,

$\bullet$ $Der\ A=\left\{
\left[\begin{array}{rcc}x&0&0\\y&2x&u\\z&0&v\end{array}\right]\ |\
x,y,z,u,v\in \mathbb{R}\right\}$,

$\bullet$ $Aut\ A=\left\{
\left[\begin{array}{rcc}x&0&0\\y&x^{2}&u\\z&0&v\end{array}\right]\
|\ x,y,z,u,v\in \mathbb{R},\ xv\neq 0\}\right.$,

$\bullet$ the partition $\mathcal{P}_{A}$ of $\mathbb{R}^{3}$,
defined by the lattice of subalgebras of $A$, consists of:

\vspace{1mm} \hspace{3mm} $\diamond$  the singletons covering the
plane $x^{2}Ox^{3}$,

\vspace{1mm} \hspace{3mm} $\diamond$ the half-planes delimited by
axis $Ox^{2}$ on each plane passing through axis $Ox^{2}$,

\vspace{1mm}$\bullet$ the partition $\mathcal{P}_{A}$ of $A$
induces a partition on the set of integral curves of the
associated homogeneous quadratic differential system (HQDS)
consisting of:

\vspace{1mm} \hspace{5mm} $\diamond$  the singletons consisting of
singular solutions that cover the plane $x^{2}Ox^{3}$,

\vspace{1mm} \hspace{5mm} $\diamond$ the integral curves contained
in the half-planes delimited by axis $Ox^{2}$ on each plane
passing through axis $Ox^{2}$.

Note that each nonsingular integral curve has both torsion and
curvature tensors zero; indeed, each of them is lying on a line of
the form $x^{1}=k_{1},\ x^{3}=k_{3}$.

\vspace{5mm}\noindent \textbf{A33)} \emph{Properties of algebras
$A=A_{7}(1,\beta)$ (of type $(ii_{2})$)} when $\beta\notin\{0,
1\}$

\vspace{3mm}$\bullet$ the algebras $A_{6}(1,\beta)\
(\beta\notin\{0, 1\})$\ and\ $A_{6}(1,\beta')\ (\beta'\notin\{0,
1\})$ are isomorphic if and only if\ $\beta=\beta '$,

$\bullet$ $Ann\ A=\{0\},\ \mathcal{N}(A)=\mathbb{R}e_{2}\cup
\mathbb{R}e_{3},\ \mathcal{I}(A)=\emptyset$,

$\bullet$ 1-dimensional subalgebras: $\mathbb{R}e_{2},\
\mathbb{R}e_{3}$,

$\bullet$ 2-dimensional subalgebras:
$Span_{\mathbb{R}}\{e_{1},e_{2}\},\
Span_{\mathbb{R}}\{e_{2},e_{3}\}$,

$\bullet$ ideals: $\mathbb{R}e_{2},\
Span_{\mathbb{R}}\{e_{1},e_{2}\}$,

$\bullet$ $A^{2}=\ Span_{\mathbb{R}}\{e_{1},e_{2}\}$; $A$ is
solvable; $A/A^{2}$ is a 1-dimensional null algebra,

$\bullet$ $Der\ A=\left\{
\left[\begin{array}{rcc}x&0&0\\0&2x&0\\0&0&0\end{array}\right]\ |\
x\in \mathbb{R}\right\}$,

$\bullet$ $Aut\ A=\left\{
\left[\begin{array}{rcc}x&0&0\\0&x^{2}&0\\0&0&1\end{array}\right]\
|\ x\in \mathbb{R}^{\ast}\right\}\cong \mathbb{R}^{\ast}(\cdot)$,

$\bullet$ the partition $\mathcal{P}_{A}$ of $\mathbb{R}^{3}$,
defined by the lattice of subalgebras of $A$, consists of:

\vspace{1mm} \hspace{3mm} $\diamond$  the singletons covering the
axes $Ox^{2}$ and $Ox^{3}$,

\vspace{1mm} \hspace{3mm} $\diamond$ the half-planes delimited by
axis $Ox^{2}$ on $x^{1}Ox^{2}$,

\vspace{1mm} \hspace{3mm} $\diamond$ the quarters of plane
delimited by axis $Ox^{2}$ and $Ox^{3}$ on $x^{2}Ox^{3}$,

\vspace{1mm} \hspace{5mm} $\diamond$ the quarters of space
delimited by planes $x^{1}Ox^{2}$ and $x^{2}Ox^{3}$,

\vspace{1mm}$\bullet$ the partition $\mathcal{P}_{A}$ of $A$
induces a partition on the set of integral curves of the
associated homogeneous quadratic differential system (HQDS)
consisting of:

\vspace{1mm} \hspace{5mm} $\diamond$  the singletons consisting of
singular solutions that cover the axes $Ox^{2}$ and $Ox^{3}$,

\vspace{1mm} \hspace{5mm} $\diamond$ the integral curves contained
in the half-planes delimited by axis $Ox^{2}$ on $x^{1}Ox^{2}$,

\vspace{1mm} \hspace{5mm} $\diamond$ the integral curves contained
in the quarters of plane delimited by axis $Ox^{2}$ and $Ox^{3}$
on $x^{2}Ox^{3}$,

\vspace{1mm} \hspace{5mm} $\diamond$ the integral curves contained
in the quarters of space delimited by planes $x^{1}Ox^{2}$ and
$x^{2}Ox^{3}$.

Note that each nonsingular integral curve are torsion-free.

\vspace{5mm}\noindent \textbf{A34)} \emph{Properties of algebras
$A=A_{7}(1,1)$ (of type $(ii_{2})$)}

\vspace{2mm} $\bullet$ $Ann\ A=\{0\},\
\mathcal{N}(A)=\mathbb{R}e_{2}\cup \mathbb{R}e_{3},\
\mathcal{I}(A)=\emptyset$,

$\bullet$ 1-dimensional subalgebras: $\mathbb{R}e_{2},\
\mathbb{R}e_{3}$,

$\bullet$ 2-dimensional subalgebras:
$Span_{\mathbb{R}}\{e_{1},e_{2}\},\
Span_{\mathbb{R}}\{e_{2},e_{3}\}$,

$\bullet$ ideals: $\mathbb{R}e_{2},\
Span_{\mathbb{R}}\{e_{1},e_{2}\}$,

$\bullet$ $A^{2}=\ Span_{\mathbb{R}}\{e_{1},e_{2}\}$; $A$ is
solvable; $A/A^{2}$ is a 1-dimensional null algebra,

$\bullet$ $Der\ A=\left\{
\left[\begin{array}{rcc}x&0&0\\y&2x&0\\0&0&0\end{array}\right]\ |\
x,y\in \mathbb{R}\right\}$,

$\bullet$ $Aut\ A=\left\{
\left[\begin{array}{rcc}x&0&0\\y&x^{2}&0\\0&0&1\end{array}\right]\
|\ x,y\in \mathbb{R},\ x\neq 0\}\right.$,

$\bullet$ the partition $\mathcal{P}_{A}$ of $\mathbb{R}^{3}$,
defined by the lattice of subalgebras of $A$, consists of:

\vspace{1mm} \hspace{3mm} $\diamond$  the singletons covering the
axes $Ox^{2}$ and $Ox^{3}$,

\vspace{1mm} \hspace{3mm} $\diamond$ the half-planes delimited by
axis $Ox^{2}$ on $x^{1}Ox^{2}$,

\vspace{1mm} \hspace{3mm} $\diamond$ the quarters of plane
delimited by axis $Ox^{2}$ and $Ox^{3}$ on $x^{2}Ox^{3}$,

\vspace{1mm} \hspace{5mm} $\diamond$ the quarters of space
delimited by planes $x^{1}Ox^{2}$ and $x^{2}Ox^{3}$,

\vspace{1mm}$\bullet$ the partition $\mathcal{P}_{A}$ of $A$
induces a partition on the set of integral curves of the
associated homogeneous quadratic differential system (HQDS)
consisting of:

\vspace{1mm} \hspace{5mm} $\diamond$  the singletons consisting of
singular solutions that cover the axes $Ox^{2}$ and $Ox^{3}$,

\vspace{1mm} \hspace{5mm} $\diamond$ the integral curves contained
in the half-planes delimited by axis $Ox^{2}$ on $x^{1}Ox^{2}$,

\vspace{1mm} \hspace{5mm} $\diamond$ the integral curves contained
in the quarters of plane delimited by axis $Ox^{2}$ and $Ox^{3}$
on $x^{2}Ox^{3}$,

\vspace{1mm} \hspace{5mm} $\diamond$ the integral curves contained
in the quarters of space delimited by planes $x^{1}Ox^{2}$ and
$x^{2}Ox^{3}$.

Note that each nonsingular integral curve are torsion-free.

\vspace{5mm}\noindent \textbf{A35)} \emph{Properties of algebras
$A=A_{7}(1,0)$ (of type $(ii_{2})$)}

\vspace{2mm} $\bullet$ $Ann\ A=\mathbb{R}e_{2},\
\mathcal{N}(A)=Span_{\mathbb{R}}\{e_{2},e_{3}\},\
\mathcal{I}(A)=\emptyset$,

$\bullet$ 1-dimensional subalgebras: $\mathbb{R}u$ for $u\in
\mathcal{N}(A)$,

$\bullet$ 2-dimensional subalgebras:
$Span_{\mathbb{R}}\{e_{1},e_{2}\},\
Span_{\mathbb{R}}\{e_{2},e_{3}\}$,

$\bullet$ ideals: $\mathbb{R}e_{2},\
Span_{\mathbb{R}}\{e_{1},e_{2}\}$,

$\bullet$ $A^{2}=\ Span_{\mathbb{R}}\{e_{1},e_{2}\}$; $A$ is
solvable; $A/A^{2}$ is the null 1-dimensional algebra,

$\bullet$ $Der\ A=\left\{
\left[\begin{array}{rcc}x&0&0\\0&2x&y\\0&0&0\end{array}\right]\ |\
x,y\in \mathbb{R}\right\}$,

$\bullet$ $Aut\ A=\left\{
\left[\begin{array}{rcc}x&0&0\\0&x^{2}&y\\0&0&1\end{array}\right]\
|\ x,y\in \mathbb{R},\ x\neq 0\}\right.$,

$\bullet$ the partition $\mathcal{P}_{A}$ of $\mathbb{R}^{3}$,
defined by the lattice of subalgebras of $A$, consists of:

\vspace{1mm} \hspace{3mm} $\diamond$  the singletons covering the
plane $x^{2}Ox^{3}$,

\vspace{1mm} \hspace{3mm} $\diamond$ the half-planes delimited by
axis $Ox^{2}$ on $x^{1}Ox^{2}$,

\vspace{1mm} \hspace{5mm} $\diamond$ the quarters of space
delimited by planes $x^{1}Ox^{2}$ and $x^{2}Ox^{3}$,

\vspace{1mm}$\bullet$ the partition $\mathcal{P}_{A}$ of $A$
induces a partition on the set of integral curves of the
associated homogeneous quadratic differential system (HQDS)
consisting of:

\vspace{1mm} \hspace{5mm} $\diamond$  the singletons consisting of
singular solutions that cover the plane $x^{2}Ox^{3}$,

\vspace{1mm} \hspace{5mm} $\diamond$ the integral curves contained
in the half-planes delimited by axis $Ox^{2}$ on $x^{1}Ox^{2}$,

\vspace{1mm} \hspace{5mm} $\diamond$ the integral curves contained
in the quarters of space delimited by planes $x^{1}Ox^{2}$ and
$x^{2}Ox^{3}$.

Note that each nonsingular integral curve are torsion-free.

\vspace{5mm}\noindent \emph{Properties of algebras $(ii_{3})$}

\begin{proposition} Each algebra of type $A_{7}(\alpha,1)$ for a given $\alpha\neq 0$
is isomorphic to the algebra of type $A_{7}(1,\beta)$ for $\beta=
\frac{1}{\alpha}$.\end{proposition}

\vspace{3mm} \noindent\emph{Subcase (iii)\ b=j=0}

\vspace{3mm} There exists a basis
$\mathcal{B}=(e_{1},e_{2},e_{3})$ such that the multiplication
table of algebra has one of the following two multiplication
tables:

$$\begin{array}{llll}
  (iii_{1})& \hspace{5mm} e_{1}^{2}=0  &\hspace{5mm} e_{2}^{2}=0&\hspace{5mm} e_{3}^{2}=0 \\
  \hspace{5mm} & \hspace{5mm} e_{1}e_{2}=0&\hspace{5mm}  e_{1}e_{3}=e_{1} &\hspace{5mm}  e_{2}e_{3}=\beta e_{2}\\
  \\
(iii_{2})& \hspace{5mm} e_{1}^{2}=0  &\hspace{5mm} e_{2}^{2}=0&\hspace{5mm} e_{3}^{2}=0 \\
  \hspace{5mm} & \hspace{5mm} e_{1}e_{2}=0&\hspace{5mm}  e_{1}e_{3}=\alpha e_{1} &\hspace{5mm}  e_{2}e_{3}=e_{2}
\end{array}$$

\vspace{5mm}\noindent \emph{Properties of algebras $(iii_{1})$}

\vspace{3mm} There exists a basis
$\mathcal{B}=(e_{1},e_{2},e_{3})$ such that the multiplication
table of algebra becomes:
$$\begin{array}{llll}
  (ii_{1})& \hspace{5mm} e_{1}^{2}=0  &\hspace{5mm} e_{2}^{2}=0&\hspace{5mm} e_{3}^{2}=0 \\
  \hspace{5mm} & \hspace{5mm} e_{1}e_{2}=0&\hspace{5mm}  e_{1}e_{3}=e_{1} &\hspace{5mm}  e_{2}e_{3}=\beta e_{2}
\end{array}$$
Consequently, this algebra is isomorphic to algebra \textbf{A22}).

Moreover, the following result is true.
\begin{proposition} Each algebra of type $(iii_{2})$ for a given $\alpha\neq 0$
is isomorphic to the algebra of type $(iii_{1})$ for $\beta=
\alpha$.\end{proposition}

\vspace{3mm}\textbf{1) Case} $Spec\ D=(1,\omega,0)$ for
$\omega\notin \{-1,0,\frac{1}{2},1,2\}$

\vspace{3mm}Algebra $A(\cdot)$ has the next multiplication table:
$$\begin{array}{llll}
  \textbf{Table T}\hspace{6mm} & \hspace{5mm} e_{1}^{2}=0  &\hspace{5mm} e_{2}^{2}=0&\hspace{5mm} e_{3}^{2}=je_{3} \\
   &\hspace{5mm}  e_{1}e_{2}=0  &\hspace{5mm}  e_{1}e_{3}=pe_{1} &\hspace{5mm}  e_{2}e_{3}=te_{2}
\end{array}$$
with $j,p,t\in \mathbb{R}$. Then, in basis $(e_{1}, e_{2},
\frac{1}{j}e_{3})$, the multiplication table of algebra becomes:
$$\begin{array}{llll}
  \textbf{Table T}\hspace{6mm} & \hspace{5mm} e_{1}^{2}=0  &\hspace{5mm} e_{2}^{2}=0&\hspace{5mm} e_{3}^{2}=e_{3} \\
   &\hspace{5mm}  e_{1}e_{2}=0  &\hspace{5mm}  e_{1}e_{3}=\alpha e_{1} &\hspace{5mm}  e_{2}e_{3}=\beta e_{2}
\end{array}$$
with $\alpha,\beta\in \mathbb{R}$. Consequently, each such algebra
is isomorphic to algebras of type $\textbf{A14}$).

\vspace{3mm}Theorem \ref{t1} stated that any algebra in class
$\textbf{Ai}$ is not isomorphic to any algebra of class
$\textbf{Aj}$ for $i,j\in \{1,2,...,7\}$. This result can be now
improved by a straight comparison of lists of properties for
isomorphism classes $\textbf{\textbf{Ai}}$ for $i\in
\{1,2,...,35\}$.

\begin{theorem}\label{t2} Each algebra of type $\textbf{Ai}$ is not isomorphic to any algebra of type
$\textbf{Aj}$ for $i,j\in\{1,2,...,35\}$ and $i\neq
j$.\end{theorem} This result induces a classification result for
the corresponding HQDSs up to an affine-equivalence.

\begin{theorem}\label{t3}For any nontrivial HQDS on $\mathbb{R}^{3}$, having at least a semisimple derivation with 1-dimensional kernel, there
exists a center-affinity such that it turns into one of the
following 35 HQDSs:

$$\begin{array}{ll}{\hspace{2mm}1^{\circ}.\left\{\begin{array}{l}\frac{\ds dx^{1}}{\ds dt}=0\\ \frac{\ds dx^{2}}{\ds dt}=0\\
\frac{\ds dx^{3}}{\ds dt}=2x^{1}x^{2}+(x^{3})^{2},
\end{array}\right.}&{\hspace{2mm}2^{\circ}.\left\{\begin{array}{l}\frac{\ds dx^{1}}{\ds dt}=0\\
\frac{\ds dx^{2}}{\ds dt}=x^{2}x^{3}\\ \frac{\ds dx^{3}}{\ds
dt}=2x^{1}x^{2}+(x^{3})^{2},
\end{array}\right.} \\
\\
{\hspace{2mm}3^{\circ}.\left\{ \begin{array}{l}\frac{\ds dx^{1}}{\ds dt}=x^{1}x^{3}\\
\frac{\ds dx^{2}}{\ds dt}=x^{2}x^{3}\\  \frac{\ds dx^{3}}{\ds
dt}=2x^{1}x^{2}+(x^{3})^{2},\end{array}\right.
}&{\hspace{2mm}4^{\circ}.\left\{\begin{array}{l}\frac{\ds
dx^{1}}{\ds
dt}=2\alpha x^{1}x^{3}\\ \frac{\ds dx^{2}}{\ds dt}=0\\
\frac{\ds dx^{3}}{\ds dt}=2x^{1}x^{2}+(x^{3})^{2}, \end{array}\right.}\\
\hspace{12mm}&\hspace{12mm}(\alpha\notin\{0,\frac{1}{2}\})
\\
\\
{\hspace{2mm}5^{\circ}.\left\{
\begin{array}{l}\frac{\ds dx^{1}}{\ds dt}=2\alpha x^{1}x^{3}\\ \frac{\ds
dx^{2}}{\ds dt}=x^{2}x^{3}\\  \frac{\ds dx^{3}}{\ds
dt}=2x^{1}x^{3}+(x^{3})^{2},\end{array}\right. }&
{\hspace{2mm}6^{\circ}.\left\{\begin{array}{l}\frac{\ds dx^{1}}{\ds dt}=2\alpha x^{1}x^{3}\\
\frac{\ds dx^{2}}{\ds dt}=2\alpha x^{2}x^{3}\\ \frac{\ds
dx^{3}}{\ds dt}=2x^{1}x^{2}+(x^{3})^{2},
\end{array}\right.}\\
\hspace{12mm}(\alpha\notin\{0,\frac{1}{2}\}) &\hspace{12mm}(\alpha\notin\{0,\frac{1}{2}\})\\
\\
{\hspace{2mm}7^{\circ}.\left\{\begin{array}{l}\frac{\ds
dx^{1}}{\ds
dt}=2\alpha x^{1}x^{3}\\ \frac{\ds dx^{2}}{\ds dt}=2\beta x^{2}x^{3}\\
\frac{\ds dx^{3}}{\ds dt}=2x^{1}x^{2}+(x^{3})^{2},
\end{array}\right.} & {\hspace{2mm}8^{\circ}.\left\{
\begin{array}{l}\frac{\ds dx^{1}}{\ds dt}=0\\ \frac{\ds
dx^{2}}{\ds dt}=0\\  \frac{\ds dx^{3}}{\ds
dt}=(x^{3})^{2},\end{array}\right. }\\
\hspace{12mm}(\alpha,\beta\notin\{0,\frac{1}{2}\},\ \alpha<\beta)&
\end{array}$$

$$\begin{array}{ll}{\hspace{2mm}9^{\circ}.\left\{\begin{array}{l}\frac{\ds
dx^{1}}{\ds
dt}=0\\ \frac{\ds dx^{2}}{\ds dt}= x^{2}x^{3}\\
\frac{\ds dx^{3}}{\ds dt}=(x^{3})^{2},
\end{array}\right.}& {10^{\circ}.\left\{
\begin{array}{l}\frac{\ds dx^{1}}{\ds dt}=x^{1}x^{3}\\ \frac{\ds
dx^{2}}{\ds dt}=x^{2}x^{3}\\  \frac{\ds dx^{3}}{\ds
dt}=(x^{3})^{2},\end{array}\right. }\\
\\

{11^{\circ}.\left\{\begin{array}{l}\frac{\ds dx^{1}}{\ds dt}=0\\
\frac{\ds dx^{2}}{\ds dt}=x^{2}x^{3}\\ \frac{\ds dx^{3}}{\ds
dt}=(x^{3})^{2},
\end{array}\right.}
&{12^{\circ}.\left\{\begin{array}{l}\frac{\ds dx^{1}}{\ds dt}=x^{1}x^{3}\\
\frac{\ds dx^{2}}{\ds dt}=2\beta x^{2}x^{3}\\ \frac{\ds
dx^{3}}{\ds dt}=(x^{3})^{2},
\end{array}\right.}\\
&\hspace{12mm}(\beta\notin\{0,\frac{1}{2}\})
\\
\\
{13^{\circ}.\left\{\begin{array}{l}\frac{\ds dx^{1}}{\ds dt}=2\alpha x^{1}x^{3}\\
\frac{\ds dx^{2}}{\ds dt}=2\alpha x^{2}x^{3}\\ \frac{\ds
dx^{3}}{\ds dt}=(x^{3})^{2},
\end{array}\right.}&{14^{\circ}.\left\{\begin{array}{l}\frac{\ds dx^{1}}{\ds dt}=2\alpha x^{1}x^{3}\\
\frac{\ds dx^{2}}{\ds dt}=2\beta x^{2}x^{3}\\ \frac{\ds
dx^{3}}{\ds dt}=(x^{3})^{2},
\end{array}\right.}\\
\hspace{12mm}(\alpha\notin\{0,\frac{1}{2}\})&\hspace{12mm}(\alpha,
\beta\notin\{0,\frac{1}{2}\},\ \alpha< \beta)
\\
\\
{15^{\circ}.\left\{\begin{array}{l}\frac{\ds dx^{1}}{\ds
dt}=2x^{2}x^{3}\\ \frac{\ds dx^{2}}{\ds dt}=0\\
\frac{\ds dx^{3}}{\ds dt}=0,
\end{array}\right.}&{16^{\circ}.\left\{\begin{array}{l}\frac{\ds dx^{1}}{\ds dt}=0\\
\frac{\ds dx^{2}}{\ds dt}=2 x^{2}x^{3}\\ \frac{\ds dx^{3}}{\ds
dt}=2 x^{1}x^{2},
\end{array}\right.}\\
\\
{17^{\circ}.\left\{\begin{array}{l}\frac{\ds dx^{1}}{\ds dt}=2x^{1}x^{3}\\
\frac{\ds dx^{2}}{\ds dt}=2 x^{2}x^{3}\\ \frac{\ds dx^{3}}{\ds
dt}=2 x^{1}x^{2},
\end{array}\right.}&{18^{\circ}.\left\{\begin{array}{l}\frac{\ds dx^{1}}{\ds dt}=2 x^{1}x^{3}\\
\frac{\ds dx^{2}}{\ds dt}=2\beta x^{2}x^{3}\\
\frac{\ds dx^{3}}{\ds dt}=2 x^{1}x^{2},
\end{array}\right.}\\
&\hspace{12mm}(\beta > 1)\\
\\
{19^{\circ}.\left\{\begin{array}{l}\frac{\ds dx^{1}}{\ds dt}=0\\
\frac{\ds dx^{2}}{\ds dt}=2 x^{2}x^{3}\\
\frac{\ds dx^{3}}{\ds dt}=0,
\end{array}\right.}&{20^{\circ}.\left\{\begin{array}{l}\frac{\ds dx^{1}}{\ds dt}=2x^{1}x^{3}\\
\frac{\ds dx^{2}}{\ds dt}=2 x^{2}x^{3}\\
\frac{\ds dx^{3}}{\ds dt}=0,\end{array}\right.}
\\
\\
{21^{\circ}.\left\{\begin{array}{l}\frac{\ds dx^{1}}{\ds dt}=2x^{1}x^{3}\\
\frac{\ds dx^{2}}{\ds dt}=2\beta x^{2}x^{3}\\
\frac{\ds dx^{3}}{\ds dt}=0,
\end{array}\right.}&{22^{\circ}.\left\{\begin{array}{l}\frac{\ds dx^{1}}{\ds dt}=2x^{2}x^{3}\\
\frac{\ds dx^{2}}{\ds dt}=0\\
\frac{\ds dx^{3}}{\ds dt}=(x^{3})^{2},\end{array}\right.}
\\
\hspace{12mm}(\beta \notin \{0, 1\})&\\
\\
{23^{\circ}.\left\{\begin{array}{l}\frac{\ds dx^{1}}{\ds dt}=2ax^{1}x^{2}+2bx^{2}x^{3}\\
\frac{\ds dx^{2}}{\ds dt}=-2b x^{1}x^{3}+2a x^{2}x^{3}\\
\frac{\ds dx^{3}}{\ds dt}=(x^{3})^{2},
\end{array}\right.}&{24^{\circ}.\left\{\begin{array}{l}\frac{\ds dx^{1}}{\ds dt}=(x^{2})^{2}\\
\frac{\ds dx^{2}}{\ds dt}=2\alpha x^{2}x^{3}\\
\frac{\ds dx^{3}}{\ds dt}=(x^{3})^{2},\end{array}\right.}\\
\hspace{12mm}(b>0)&\hspace{12mm}(\alpha\notin\{0,\frac{1}{2}\})
\end{array}$$

$$\begin{array}{ll}
{25^{\circ}.\left\{\begin{array}{l}\frac{\ds dx^{1}}{\ds dt}=(x^{2})^{2}\\
\frac{\ds dx^{2}}{\ds dt}=0\\
\frac{\ds dx^{3}}{\ds dt}=(x^{3})^{2},
\end{array}\right.}& {26^{\circ}.\left\{\begin{array}{l}\frac{\ds dx^{1}}{\ds dt}=(x^{2})^{2}\\
\frac{\ds dx^{2}}{\ds dt}=x^{2}x^{3}\\
\frac{\ds dx^{3}}{\ds
dt}=(x^{3})^{2},\end{array}\right.}\\
\\

{27^{\circ}.\left\{\begin{array}{l}\frac{\ds dx^{1}}{\ds dt}=2\alpha x^{1}x^{3}+(x^{2})^{2}\\
\frac{\ds dx^{2}}{\ds dt}=2\beta x^{2}x^{3}\\
\frac{\ds dx^{3}}{\ds dt}=(x^{3})^{2},
\end{array}\right.}&{28^{\circ}.\left\{\begin{array}{l}\frac{\ds dx^{1}}{\ds dt}=2\alpha x^{1}x^{3}+(x^{2})^{2}\\
\frac{\ds dx^{2}}{\ds dt}=2\alpha x^{2}x^{3}\\
\frac{\ds dx^{3}}{\ds dt}=(x^{3})^{2},\end{array}\right.}
\\
\hspace{12mm}(\alpha,\ \beta \notin \{0, \frac{1}{2}\},\
\alpha\neq \beta)&\hspace{12mm}(\alpha \notin \{0, \frac{1}{2}\})
\\
\\

{29^{\circ}.\left\{\begin{array}{l}\frac{\ds dx^{1}}{\ds dt}=2\alpha x^{1}x^{3}+(x^{2})^{2}\\
\frac{\ds dx^{2}}{\ds dt}=x^{2}x^{3}\\
\frac{\ds dx^{3}}{\ds dt}=(x^{3})^{2},
\end{array}\right.}&{30^{\circ}.\left\{\begin{array}{l}\frac{\ds dx^{1}}{\ds dt}= x^{1}x^{3}+(x^{2})^{2}\\
\frac{\ds dx^{2}}{\ds dt}= x^{2}x^{3}\\
\frac{\ds dx^{3}}{\ds dt}=(x^{3})^{2},\end{array}\right.}
\\
\hspace{12mm}(\alpha \notin \{0, \frac{1}{2}\})&\\
\\
\\
{31^{\circ}.\left\{\begin{array}{l}\frac{\ds dx^{1}}{\ds dt}=x^{1}x^{3}+(x^{2})^{2}\\
\frac{\ds dx^{2}}{\ds dt}=2\beta x^{2}x^{3}\\
\frac{\ds dx^{3}}{\ds dt}=(x^{3})^{2},
\end{array}\right.}&{32^{\circ}.\left\{\begin{array}{l}\frac{\ds dx^{1}}{\ds dt}=(x^{1})^{2}\\
\frac{\ds dx^{2}}{\ds dt}=0\\
\frac{\ds dx^{3}}{\ds dt}=0,\end{array}\right.}
\\
\hspace{12mm}(\beta \notin \{0, \frac{1}{2}\})&
\\
\\
{33^{\circ}.\left\{\begin{array}{l}\frac{\ds dx^{1}}{\ds dt}=2x^{1}x^{3}\\
\frac{\ds dx^{2}}{\ds dt}=(x^{1})^{2}+2\beta x^{2}x^{3}\\
\frac{\ds dx^{3}}{\ds dt}=0,
\end{array}\right.}&{34^{\circ}.\left\{\begin{array}{l}\frac{\ds dx^{1}}{\ds dt}=2 x^{1}x^{3}\\
\frac{\ds dx^{2}}{\ds dt}=(x^{1})^{2}+2x^{2}x^{3}\\
\frac{\ds dx^{3}}{\ds dt}=0\end{array}\right.},\\
\hspace{12mm}(\beta \notin \{0, 1\})&
\\
\\
{35^{\circ}.\left\{\begin{array}{l}\frac{\ds dx^{1}}{\ds dt}=2x^{1}x^{3}\\
\frac{\ds dx^{2}}{\ds dt}=(x^{1})^{2}\\
\frac{\ds dx^{3}}{\ds dt}=0.
\end{array}\right.}&\

\end{array}$$
\end{theorem}

\section{Conclusions}

The existence of a semisimple derivation for a HQDS is a strong
constraint implying the existence of a center-affine equivalent
system having a lot of coefficients either 0 or small integers in
suitable bases. The classification result of this family of HQDSs
up to an affine equivalence, just exhibited in Theorem \ref{t3},
is a corollary of the classification of corresponding commutative
algebras up to an isomorphism. There were identified 13 classes of
algebras $A$ having $Der\ A=\mathbb{R}D$ while their automorphism
groups are isomorphic either with $\mathbb{R}^{\ast}(\cdot)$ or to
$\mathbb{R}^{\ast}\times \{-1,1\}$. The other 21 classes of
algebras have larger derivations algebras and automorphism groups.
Note that, for 3 algebras the automorphism group is
$\mathbb{R}^{\ast}(\cdot)\times \mathbb{R}^{\ast}(\cdot)$, for 3
algebras the automorphism group is $GL(2,\mathbb{R})$, one algebra
has $Aff(2,\mathbb{R})$ as its automorphism group and one algebra
has $\mathbb{C}^{\ast}$ as its automorphism group; the other 14
algebras have more complex automorphism groups. In order to decide
on the mutually non-isomorphism of pairs of algebras in different
classes $\textbf{Ai}$ for $i\in\{1,2,...,35\}$ we have supplied
the lists of main properties assigned to each such a class. Note
that the algebras having the same list of main properties were
collected together by means of a label consisting in one or two
parameters running over a specific range. In fact, this
classification seems to be the finest under the isomorphism
criterion. In these lists were also included the information
concerning the partitions of ground space $A$ as well as of the
set of all integral curves of the associated systems. These
partitions are invariant under the action of automorphism groups
so that they can be used as tests of correctness of results. On
the other hand we hope that these partitions will give important
information about the stability of steady state solutions.
Finally, let us remark that some algebraic properties, like the
solvability or the existence of special ideals of algebras,
attract important geometric properties of nonsingular integral
curves of the corresponding HQDS like the vanishing of torsion or
curvature tensor.

\end{document}